\documentclass[a4paper]{article}

\usepackage[T1]{fontenc}
\usepackage[utf8]{inputenc}
\usepackage[english]{babel}

\usepackage{amssymb,amsmath,amsfonts,amsthm}   

\usepackage{mathtools}   
\usepackage{bm}   

\usepackage{graphicx}

\usepackage{tikz}
\usetikzlibrary{matrix,arrows,decorations.pathmorphing}
\usepackage{tikz-cd}

\usepackage[mathscr]{eucal}   
\usepackage[sans]{dsfont}   

\usepackage{enumitem}   

\usepackage{afterpage}

\usepackage[nottoc]{tocbibind}   

\theoremstyle{plain} 
\newtheorem{thm}{Theorem}[section]
\newtheorem{prop}[thm]{Proposition}
\newtheorem{lmm}[thm]{Lemma}
\newtheorem{cor}[thm]{Corollary}

\theoremstyle{definition} 

\theoremstyle{remark} 
\newtheorem{rmk}[thm]{Remark}
\newtheorem{de}[thm]{Definition}

\newtheorem*{notat}{Notation and conventions}

\DeclareMathOperator{\Aut}{Aut}

\DeclareMathOperator{\coker}{coker}
\DeclareMathOperator{\car}{char}

\DeclareMathOperator{\Exp}{Exp}
\DeclareMathOperator{\Ext}{Ext}
\DeclareMathOperator{\h}{H}

\DeclareMathOperator{\Hom}{Hom}
\DeclareMathOperator{\id}{id}
\DeclareMathOperator{\ind}{Ind}
\DeclareMathOperator{\len}{length}
\DeclareMathOperator{\Log}{Log}
\DeclareMathOperator{\Mor}{Mor}

\DeclareMathOperator{\Quot}{Quot}

\DeclareMathOperator{\spec}{Spec}

\newcommand{\Z}{\mathbb{Z}}
\newcommand{\Q}{\mathbb{Q}}
\newcommand{\C}{\mathbb{C}}

\newcommand{\K}{\Bbbk}   
\newcommand{\aff}{\mathbb{A}}   
\newcommand{\pr}{\mathbb{P}}   
\newcommand{\GL}{\mathrm{GL}}   
\newcommand{\grass}{\mathds{G}}
\newcommand{\sym}{\mathbb{S}}

\renewcommand{\phi}{\varphi}

\newcommand{\mb}[1]{\mathbf{#1}}
\newcommand{\mc}[1]{\mathcal{#1}}

\newcommand{\mr}[1]{\mathrm{#1}}
\newcommand{\ms}[1]{\mathsf{#1}}

\newcommand{\iso}{\xrightarrow{{}_\thicksim}}
\newcommand{\into}{\hookrightarrow}

\newcommand{\xto}{\xrightarrow}

\newcommand{\1}{{1\!\!1}}   

\newcommand{\cls}[1]{\mkern 1.5mu\overline{\mkern-1.5mu#1\mkern-1.5mu}\mkern 1.5mu}   

\makeatletter
\newcommand*{\bigcdot}{}
\DeclareRobustCommand*{\bigcdot}{%
  \mathbin{\mathpalette\bigcdot@{}}%
}
\newcommand*{\bigcdot@scalefactor}{.5}
\newcommand*{\bigcdot@widthfactor}{1.15}
\newcommand*{\bigcdot@}[2]{%
  \sbox0{$#1\vcenter{}$}
  \sbox2{$#1\cdot\m@th$}%
  \hbox to \bigcdot@widthfactor\wd2{%
    \hfil
    \raise\ht0\hbox{%
      \scalebox{\bigcdot@scalefactor}{%
        \lower\ht0\hbox{$#1\bullet\m@th$}%
      }%
    }%
    \hfil
  }%
}
\makeatother

\newcommand{\Kab}{K_0(\ms{A})}
\newcommand{\eKab}{K_0^{\sym}(\ms{A})}
\newcommand{\gKab}{K_0(\ms{A}^{\Upsilon})}
\newcommand{\geKab}{K_0^{\sym}(\ms{A}^{\Upsilon})}

\newcommand{\Kvar}{K_0(\ms{Var})}
\newcommand{\eKvar}{K_0^\sym(\ms{Var})}
\newcommand{\gKvar}{K_0(\ms{Var}^{\Upsilon})}
\newcommand{\geKvar}{K_0^{\sym}(\ms{Var}^{\Upsilon})}

\newcommand{\rKvar}[1]{K_0(\ms{Var}_{#1})}

\newcommand{\Kmhs}{K_0(\ms{MHS})}
\newcommand{\eKmhs}{K_0^\sym(\ms{MHS})}

\newcommand{\geKmhs}{K_0^{\sym}(\ms{MHS}^{\Upsilon})}

\newcommand{\nKmhs}[1]{K_0^{\sym_{#1}}(\ms{MHS})}
\newcommand{\nKvar}[1]{K_0^{\sym_{#1}}(\ms{Var})}

\newcommand{\maff}{\mathbb{L}}
\newcommand{\caff}{\mathds{L}}
\newcommand{\hg}{\mathfrak{e}}
\newcommand{\ser}{[\![q]\!]}

\newcommand{\Addresses}{{
  \bigskip

  \noindent Massimo Bagnarol -- \texttt{mbagnarol@protonmail.com}
  
}}

\title{Betti numbers of stable map spaces to Grassmannians}

\author{Massimo Bagnarol}

\date{}

\begin{document}


\maketitle


\begin{abstract}
Let $\cls{M}_{0,n}(G(r,V), d)$ be the coarse moduli space of stable degree $d$ maps from $n$-pointed genus $0$ curves to a Grassmann variety $G(r,V)$. We provide a recursive method for the computation of the Hodge numbers and the Betti numbers of $\cls{M}_{0,n}(G(r,V), d)$ for all $n$ and $d$. Our method is a generalization of Getzler and Pandharipande's work for maps to projective spaces.
\end{abstract}

\tableofcontents

\section{Introduction}

The relevance of Kontsevich's moduli spaces of stable maps (see \cite{Kon95}) has by now become widespread, notably in connection with Gromov-Witten theory. These spaces are compactifications of spaces of morphisms from nonsingular curves to a fixed variety, which are obtained by allowing the domain curves to have nodes.

Let us briefly recall their definition in the genus $0$ complex case, which is the one we deal with in this paper. For an exhaustive and general treatment, we refer to \cite{BM96} and \cite{FP97}. Let $Y$ be a nonsingular projective $\C$-variety, and let $\beta \in H_2(Y,\Z)$ be a curve class. A stable map from an $n$-pointed genus $0$ curve to $Y$ of class $\beta$, or stable $(Y,0,n,\beta)$-map for short, is a triple $(C,x,f)$, where $C$ is a projective, connected, reduced, at worst nodal curve of arithmetic genus $0$, $x = (x_i)_{1 \le i \le n}$ is an $n$-tuple of pairwise distinct regular points of $C$, and $f: C \to Y$ is a morphism such that $f_\ast[C] = \beta$. These data are subject to the following stability condition: if $C'$ is an irreducible component of $C$ that is contracted to a point by $f$, then $C'$ has at least three special points (i.e., either marked points or nodes).

Stable $(Y,0,n,\beta)$-maps are naturally parametrized by a proper algebraic Deligne-Mumford stack $\cls{\mc{M}}_{0,n}(Y,\beta)$, with a projective coarse space $\cls{M}_{0,n}(Y,\beta)$. If $Y$ is convex, then the stack $\cls{\mc{M}}_{0,n}(Y,\beta)$ is nonsingular and $\cls{M}_{0,n}(Y,\beta)$ is a normal variety which is locally the quotient of a nonsingular variety by a finite group. Furthermore, the boundary of $\cls{M}_{0,n}(Y,\beta)$, i.e., the locus of singular domain curves, is a divisor with normal crossings, up to a finite group quotient. It can be given a combinatorial description similar to that of the boundary of the moduli space $\cls{M}_{0,n}$ of stable $n$-pointed genus $0$ curves.

Despite the relevance of stable map spaces, their topology is still not well understood. For instance, their cohomology ring has been studied only in few cases (see \cite{BO03} and \cite{Mu07}). In the case of genus $0$ stable maps to a flag variety $Y$, it is known that the cohomology groups are generated by tautological classes (see \cite{Opr06}). In fact, computing the Betti numbers of $\cls{M}_{0,n}(Y,\beta)$ is already a nontrivial problem. For $Y = \pr^r$, this problem was solved in \cite{GP06} by Getzler and Pandharipande, who computed the Hodge and Betti numbers of $\cls{M}_{0,n}(\pr^r,d)$ by determining its Hodge-Grothendieck characteristic (referred to as Serre characteristic in \cite{GP06}), i.e., a refined $\sym_n$-equivariant version of the $E$-polynomial.

The aim of this paper is to extend the method of \cite{GP06} to the case of stable maps to a finite dimensional Grassmann variety $G(r,V)$. For some special values of $n$ and $d$, there are some computations of the Poincaré polynomial of these spaces in the literature (see, for instance, \cite{Lop14}), which are performed using the Białynicki-Birula decomposition determined by a certain torus action. However, such computations are limited to $n=0$ and $d \le 3$. Therefore, in order to determine the Hodge and Betti numbers of $\cls{M}_{0,n}(G(r,V),d)$ for all $n$ and $d$, we developed an alternative method that generalizes \cite{GP06}. The outcome of the paper is a recursive algorithm through which the Hodge-Grothendieck characteristic $\hg(\cls{M}_{0,n}(G(r,V),d))$ can be computed, for any $V,r,n,d$.

The main point of our method is the same as that of \cite{GP06}. Using the combinatorial description of the boundary of $\cls{M}_{0,n}(G(r,V),d)$ via stable marked trees, we recursively reduce the computation of the Hodge-Grothendieck characteristic of $\cls{M}_{0,n}(G(r,V),d)$ to that of the open loci $M_{0,m}(G(r,V),\delta)$ (for $\delta \le d$ and $m \le n+d-\delta$) where the domain curves are nonsingular. This reduction is obtained by translating the combinatorial properties of the spaces under consideration into recursive relations in the Grothendieck rings $(\prod_{n \ge 0} \nKvar{n})\ser$ and $(\prod_{n \ge 0} \nKmhs{n})\ser$, via the composition structures on those rings.

The Hodge-Grothendieck characteristic of $M_{0,m}(G(r,V),\delta)$ is determined by that of the space of degree $\delta$ morphism from $\pr^1$ to $G(r,V)$, and by that of the configuration space $\mr{F}(\pr^1,m)$, which has already been calculated in \cite{Get95conf}. This way, the original problem boils down to the computation of $\hg(\Mor_\delta(\pr^1,G(r,V)))$. 

Finally, by means of the Quot scheme compactification of $\Mor_\delta(\pr^1,G(r,V))$ studied by Str{\o}mme \cite{Str87}, we provide a method to explicitly compute its Hodge-Grothendieck characteristic. Together with the previous results, this gives a recursive algorithm to determine $\hg(\cls{M}_{0,n}(G(r,V),d))$.

The paper is organized as follows. In Section \ref{sezsugligrring}, we describe the composition structures on the rings $(\prod_{n \ge 0} \nKvar{n})\ser$ and $(\prod_{n \ge 0} \nKmhs{n})\ser$.

These structures, together with the stratification of $\cls{M}_{0,n}(G(r,V),d)$ by strata corresponding to stable marked trees, are used in Section \ref{echenneso} to obtain suitable recursive relations in the above Grothendieck rings.

The calculation of the class of $\Mor_\delta(\pr^1,G(r,V))$ in the Grothendieck ring of varieties, which in turn determines its Hodge-Grothendieck characteristic, is the subject of Section \ref{capmor}.

Section \ref{esempiacaso} contains some exemplifying applications of our algorithm.

\begin{notat}
For all categories $\ms{C}$ and $\ms{D}$ such that $\ms{C}$ is small, $[\ms{C},\ms{D}]$ denotes the category of functors from $\ms{C}$ to $\ms{D}$, with natural transformations of such functors as morphisms.

For all $n \in \Z_{\ge 0}$, $\sym_n$ denotes the symmetric group on $n$ elements.

For any field $\K$, by a $\K$-variety we mean a reduced separated scheme of finite type over $\K$, not necessarily irreducible. Unless otherwise stated, the field $\C$ of complex numbers is assumed to be the ground field.

In order to distinguish between set-theoretic and scheme-theoretic disjoint union, the first is denoted by $\sqcup$, while the latter is denoted by $\amalg$.

For any Cartesian diagram
\[
\begin{tikzcd}
X \times_Y T \arrow{r} \arrow{d} & X \arrow{d} \\
T \arrow{r} & Y
\end{tikzcd}
\] 
and any morphism $\psi: \mc{F} \to \mc{G}$ of sheaves on $X$, the pullback of $\psi$ to $X \times_Y T$ is denoted by $\psi_T: \mc{F}_T \to \mc{G}_T$. In particular, if $\spec(\kappa(y)) \to Y$ is the morphism determined by a point $y \in Y$, then the pullback of $\psi$ to $X_y$ is denoted by $\psi_y: \mc{F}_y \to \mc{G}_y$.
\end{notat}

\subsubsection*{Acknowledgements}

The paper is based on my Ph.D.~thesis \cite{Bag19}, which was written while I was a student at SISSA (Internation School for Advanced Studies) in Trieste. I am deeply indebted to Fabio Perroni and Barbara Fantechi, my thesis advisors. I would also like to thank Orsola Tommasi and Dan Petersen for their helpful comments.


\section{Composition structures on Grothendieck rings} \label{sezsugligrring}

In this section, we study the algebraic setting where our calculations will take place, namely the Grothendieck rings of varieties and mixed Hodge
structures with symmetric group actions. A relevant role is played by their composition structures, which are essential to translate the combinatorial properties of stable map spaces into relations in these rings, as we will see in \S\ref{echenneso}.

Our exposition is mainly based on that of \cite{GP06}, which we slightly extend in order to provide a suitable setting for studying the Betti numbers of stable map spaces to any nonsingular, projective, convex variety. For further details and proofs, the interested reader is invited to consult the author's Ph.D.\ thesis \cite[\S 1]{Bag19}.

\subsection{Grothendieck groups of varieties} \label{Ggrpvar}

Let $\K$ be a field. For simplicity, we assume that $\K$ is algebraically closed.

\begin{de}
The Grothendieck group $\rKvar{\K}$ of $\K$-varieties is the abelian group generated by the isomorphism classes of $\K$-varieties subject to the relations
\[
[X] = [Y] + [X \setminus Y] \,,
\]
where $X$ is a $\K$-variety and $Y \subseteq X$ is a closed subvariety.
\end{de}

The Grothendieck group $\rKvar{\K}$ admits different presentations, especially if $\car(\K) = 0$ (see \cite{Bit04}). The presentation we will mainly use is the following.

\begin{prop}
The Grothendieck group of $\K$-varieties can be alternatively presented as the abelian group generated by the isomorphism classes of quasi-projective $\K$-varieties subject to the relations $[X] = [Y] + [X \setminus Y]$, where $X$ is a quasi-projective $\K$-variety and $Y \subseteq X$ is a closed subvariety.
\end{prop}

The group $\rKvar{\K}$ is actually a commutative ring, whose product is defined on generators by
\[
[X] [Y] := [X \times_\K Y] \,.
\]

The class of the affine line $\aff^1_{\K}$ in $\rKvar{\K}$ is denoted by $\maff$. As proved in \cite{Bor18}, $\maff$ is a zero divisor in $\rKvar{\K}$ when $\car(\K)=0$.

Let us recall the following well-known properties of $\rKvar{\K}$, which will be used in the next sections.

\begin{prop} \label{grringdecomp}
Let $X$ be a $\K$-variety. Assume that there is a decomposition $X = Y_1 \sqcup \dots \sqcup Y_m$, where each $Y_i$ is a locally closed subvariety of $X$. Then
\[
[X] = [Y_1] + \dots + [Y_m]
\]
in $\rKvar{\K}$.
\end{prop}

\begin{prop} \label{solitarobatriviale}
Let $f: X \to Y$ be a morphism of $\K$-varieties. If $f$ is a locally trivial fibration in the Zariski topology, with fiber $F$, then 
\[
[X] = [Y] [F]
\]
in $\rKvar{\K}$.
\end{prop}

As a consequence of Proposition \ref{grringdecomp}, the decomposition $\pr^r_{\K} = \aff^0_{\K} \sqcup \dots \sqcup \aff^r_{\K}$ implies that $[\pr^r_{\K}] = 1 + \maff + \dots + \maff^r$.

If $G$ is a finite group, there are $G$-equivariant versions of the Grothendieck group of $\K$-varieties. In this paper, we only deal with the case where $\K$ is the field $\C$ of complex numbers and $G$ is the symmetric group $\sym_n$ on $n$ elements. Hereinafter, all varieties are understood to be over $\C$, and the notation $\ms{Var}_{\C}$ is shortened to $\ms{Var}$.

\begin{de}
The Grothendieck group $\nKvar{n}$ of $\sym_n$-varieties is the abelian group generated by the isomorphism classes of quasi-projective varieties with an action of $\sym_n$, subject to the relations $[X] = [Y] + [X \setminus Y]$ whenever $Y$ is a closed $\sym_n$-invariant subvariety of $X$.
\end{de}

When studying moduli spaces of stable maps to a nonsingular projective variety $Y$, an important role is played by the group of generalized power series $\bigl( \prod_{n \ge 0} \nKvar{n} \bigr)[\![\Upsilon]\!]$, where $\Upsilon \subseteq H_2(Y,\Z)$ is the monoid of curve classes of $Y$. This group can be itself interpreted as a suitable Grothendieck group. 

Let $\sym$ be the permutation groupoid, i.e., the category whose objects are the nonnegative integers $n \in \Z_{\ge 0}$, with $\sym(m,n)=\emptyset$ if $m \neq n$ and $\sym(n,n) = \sym_n$.

\begin{de}
A $\Upsilon$-graded $\sym$-variety is a functor from $\sym \times \Upsilon$ to the category $\ms{Var}$ of quasi-projective varieties. A morphism of $\Upsilon$-graded $\sym$-varieties is a natural transformation of such functors.
\end{de}

The quasi-projectivity assumption is made to guarantee that quotients by finite group actions exist in $\ms{Var}$.

If $X$ is a $\Upsilon$-graded $\sym$-variety, $X(n,\beta)$ denotes the image of $(n,\beta) \in \Z_{\ge 0} \times \Upsilon$ under $X$; by definition, $X(n,\beta)$ is a quasi-projective variety with an $\sym_n$-action. 

\begin{de}
The Grothendieck group $\geKvar$ of $\Upsilon$-graded $\sym$-varieties is the abelian group generated by the isomorphism classes of $\Upsilon$-graded $\sym$-varieties subject to the relations 
\[
[X] = [Y] + [X \setminus Y] \,, 
\]
where $X$ is a $\Upsilon$-graded $\sym$-variety, $Y(n,\beta)$ is a closed $\sym_n$-invariant subvariety of $X(n,\beta)$ for every $(n,\beta)$, and $(X \setminus Y)(n,\beta) := X(n,\beta) \setminus Y(n,\beta)$.
\end{de}

A $\{0\}$-graded $\sym$-variety is simply called an $\sym$-variety, and the corresponding Grothendieck group is denoted by $\eKvar$. A $\Z_{\ge 0}$-graded $\sym$-variety is called a graded $\sym$-variety. 

There are canonical isomorphisms
\[
\geKvar \cong \eKvar[\![\Upsilon]\!] \cong \Biggl( \prod_{n \ge 0} \nKvar{n} \Biggr) [\![\Upsilon]\!] \,.
\]
The presentation of this group as the Grothendieck group of $\Upsilon$-graded $\sym$-varieties endows it with some additional structures, as shown in \cite{GP06}. 

\begin{rmk}
The structures we are going to describe can be traced back to Joyal's theory of combinatorial species (see \cite{Joy81} and \cite{BLL97}). Furthermore, they are closely related with Kelly's formalism for operads of \cite{Kel05}. 
\end{rmk}

For all $n \in \Z_{\ge 0}$, let $S_n$ be the $\Upsilon$-graded $\sym$-variety defined by $S_n(i,\beta) = \emptyset$ if $(i,\beta) \neq (n,0)$ and $S_n(n,0) = \spec(\C)$ with the trivial $\sym_n$-action. The element $[S_n] \in \geKvar$ is denoted by $s_n$.

The category $[\sym \times \Upsilon, \ms{Var}]$ of $\Upsilon$-graded $\sym$-varieties is symmetric monoidal, with the product given by Day's convolution
\[
(X \boxtimes Y)(n,\beta) := \coprod_{\alpha + \gamma = \beta} \coprod_{i+j=n} \ind_{\sym_i \times \sym_j}^{\sym_n}( X(i,\alpha) \times Y(j,\gamma) )
\]
and the identity object $S_0$. Here, $\ind_{\sym_i \times \sym_j}^{\sym_n}$ is defined by the left Kan extension along the inclusion of one-object categories $\sym_i \times \sym_j \to \sym_n$. As Day's convolution product preserves coproducts in both arguments, it induces a commutative ring structure on $\geKvar$, whose product is defined on generators by
\[
[X] [Y] := [X \boxtimes Y] \,.
\]

Let $\gKvar \cong \Kvar[\![\Upsilon]\!]$ be the subring of $\geKvar$ generated by the elements $[X]$ such that $X(n,\beta) = \emptyset$ for all $n>0$. Its inclusion into $\geKvar$ induces a commutative $\gKvar$-algebra structure on $\geKvar$. Furthermore, the subalgebras
\[
F_n := \left\langle \bigcup_{i=0}^n \bigl\{ [X] \mid X(j,\beta) = \emptyset\ \text{if}\ j<n-i\ \text{or}\ m_\beta < i \bigr\} \right\rangle
\]
form a decreasing filtration of $\geKvar$, and $\geKvar$ is complete with respect to this filtration. Here, for any $\beta \neq 0$, $m_\beta$ is the maximum integer $m$ such that the set $\{ (\beta_1,\dots,\beta_m) \in \Upsilon^m \mid \beta = \sum_i \beta_i,\, \beta_i \neq 0\ \forall\, i\}$ is nonempty, while $m_0=0$.

There is another partially defined monoidal structure on $[\sym \times \Upsilon, \ms{Var}]$. Its tensor product is the composition (also known as substitution or plethysm)
\[
(X \circ Y)(n,\beta) := \coprod_{\alpha + \gamma = \beta} \coprod_{i \ge 0} \bigl( X(i,\alpha) \times Y^{\boxtimes i}(n,\gamma) \bigr)/\sym_i \,,
\]
which is defined whenever $Y(0,0) = \emptyset$, and its identity object is $S_1$.

\begin{rmk}
The composition product provides a conceptual definition of an operad with values in $\ms{Var}$: a $\ms{Var}$-valued operad $\mc{P}$ with $\mc{P}(0) = \emptyset$ is a monoid in $([\sym,\ms{Var}],\circ,S_1)$. For further details, see \cite{Kel05}.
\end{rmk}

By \cite[\S 2]{GP06}, the composition product induces a well-defined associative operation $\circ: \geKvar \times F_1 \to \geKvar$ such that
\[
[X] \circ [Y] := [X \circ Y] \,.
\] 
As the functor $- \circ Y$ preserves coproducts and commutes with Day's convolution, this operation is a $\gKvar$-algebra homomorphism in the first argument.

On the other hand, the functor $X \circ -$ does not preserve coproducts. For instance, if $X = S_n$ then 
\[
S_n \circ (Y \amalg Y') \cong \coprod_{i+j=n} (S_i \circ Y) \boxtimes (S_j \circ Y')
\]

Finally, there is an endofunctor $D$ of $[\sym \times \Upsilon, \ms{Var}]$ which is defined by
\[
(DX)(n,\beta) := X(n + 1, \beta) \,,
\]
where the $\sym_n$-action on $X(n+1,\beta)$ is the restriction of the given action of $\sym_{n+1}$ to its subgroup $\{\sigma \in \sym_{n+1} \mid \sigma(n+1)=n+1\} \cong \sym_n$.
This endofunctor preserves coproducts, and $DS_n = S_{n-1}$ for all $n \ge 1$. Moreover, $D$ satisfies both Leibniz's rule
\[
D(X \boxtimes Y) \cong (DX \boxtimes Y) \amalg (X \boxtimes DY)
\]
and the chain rule
\[
D(X \circ Y) \cong (DX \circ Y) \boxtimes DY \,.
\]
Thus, it induces a $\gKvar$-derivation $D$ on $\geKvar$ such that $D[X] := [DX]$.

Altogether, the structures we have introduced form a composition algebra structure on $\geKvar$, in the sense of \cite[Def.\ 1.1]{GP06}. For the reader's convenience, we recall the relevant definition.

\begin{de}
Let $R$ be a commutative ring, and let $A$ be a commutative algebra over $R$ which is filtered by subalgebras $A=F_0 \supseteq F_1 \supseteq \dots$. A composition operation on $A$ is an operation 
\[
\circ: A \times F_1 \to A
\]
such that $F_n \circ F_m \subseteq F_{nm}$, endowed with a derivation $D$ such that $D(F_n) \subseteq F_{n-1}$, and a sequence of elements $s_n \in F_n$ for $n \in \Z_{\ge 0}$, satisfying the following axioms:
\begin{enumerate}[label=(\roman*)]
\item for fixed $b \in F_1$, the map $- \circ b: A \to A$ is an $R$-algebra endomorphism;
\item $(a \circ b_1) \circ b_2 = a \circ (b_1 \circ b_2)$ for any $a \in A$, $b_1,b_2 \in F_1$;
\item $D(a \circ b) = (Da \circ b)(Db)$ for any $a \in A$, $b \in F_1$;
\item $s_0=1$, $s_1 \circ b = b$ for all $b \in F_1$, $a \circ s_1 = a$ for all $a \in A$, $Ds_n=s_{n-1}$, and $s_n \circ (b+b') = \sum_{i+j=n} (s_i \circ b)(s_j \circ b')$ for all $b,b' \in F_1$.
\end{enumerate}

An algebra with a composition operation is called a composition algebra. If it is also complete with respect to $F_{\bullet}$, then it is called a complete composition algebra.
\end{de}

The previous discussion can now be summarized in the following result (see \cite[Thm.\ 2.2 and \S 3]{GP06}).

\begin{thm}
The Grothendieck group of $\Upsilon$-graded $\sym$-varieties $\geKvar$ is a complete composition algebra over $\gKvar$.
\end{thm}

Since the monoid of curve classes of a finite dimensional Grassmannian is $\Z_{\ge 0}$, the only nontrivial $\Upsilon$ that will appear in the next sections is $\Upsilon = \Z_{\ge 0}$, for which we have
\[
K_0(\ms{Var}^{\Z_{\ge 0}}) \cong \Kvar\ser \qquad \text{and} \qquad K_0^{\sym}(\ms{Var}^{\Z_{\ge 0}}) \cong \eKvar\ser \,.
\] 
Nevertheless, the level of generality of the above exposition is motivated by its relevance for the study of stable maps to different targets.

\subsection{Grothendieck groups of mixed Hodge structures}

The constructions of \S\ref{Ggrpvar} are not peculiar to the category $\ms{Var}$. For instance, if $\ms{Var}$ is replaced by any abelian tensor category $(\ms{A},\otimes,\1)$ whose tensor product is exact, a composition algebra structure can be analogously defined on the Grothendieck group $\eKab$ of functors from $\sym$ to $\ms{A}$. In fact, if $\ms{A}$ is also $\Q$-linear, then $\eKab$ acquires some additional properties, which are inherited from the theory of linear representations of symmetric groups.

Our primary interest is the case where $\ms{A}$ is the category $\ms{MHS}$ of mixed Hodge structures over $\Q$. 

The relevant constructions are briefly recalled below. Our main reference is \cite[\S 5]{GP06}; for a more detailed exposition, the reader is referred to \cite{Get95conf}. Since we are interested in studying the cohomology of stable map spaces, the case of functors $\sym \times \Upsilon \to \ms{A}$ is made explicit. As in \S\ref{Ggrpvar}, $\Upsilon$ denotes the monoid of curve classes of the nonsingular projective variety $Y$ which is the target of the stable maps under consideration.

First, recall the usual definition of the Grothendieck group of an abelian category.

\begin{de} \label{Ggrpabel}
Let $\ms{C}$ be an abelian category, which is assumed to be essentially small in order to avoid set-theoretic issues. The Grothendieck group $K_0(\ms{C})$ of $\ms{C}$ is the abelian group generated by the isomorphism classes of objects of $\ms{C}$, subject to the relations
\[
[B] = [A] + [C]
\]
whenever $0 \to A \to B \to C \to 0$ is a short exact sequence in $\ms{C}$.
\end{de}

In particular, for any abelian category $\ms{A}$, Definition \ref{Ggrpabel} also applies to the functor categories $[\sym_n,\ms{A}]$, $[\sym, \ms{A}]$ and $[\sym \times \Upsilon, \ms{A}] \cong [\sym, [\Upsilon, \ms{A}]]$. The corresponding Grothendieck groups are denoted by $K_0^{\sym_n}(\ms{A})$, $\eKab$ and $\geKab$, respectively. There are canonical isomorphisms
\[
\geKab \cong \eKab[\![\Upsilon]\!] \cong \Biggl( \prod_{n \ge 0} K_0^{\sym_n}(\ms{A}) \Biggr)[\![\Upsilon]\!] \,.
\]

Now, let $(\ms{A},\otimes,\1)$ be a fixed abelian tensor category over $\Q$ (or any field of characteristic $0$), with the tensor product $\otimes$ and the identity object $\1$, such that $\ms{A}$ is essentially small and $\otimes$ is exact. The tensor product $\otimes$ induces a commutative ring structure on $\Kab$, with the identity $[\1]$.

For all $n \in \Z_{\ge 0}$, let $S_n: \sym \times \Upsilon \to \ms{A}$ be the functor defined by $S_n(i,\beta) := 0$ if $(i,\beta) \neq (n,0)$ and $S_n(n,0) := \1$ with the trivial $\sym_n$-action. The element $[S_n] \in \geKab$ is denoted by $s_n$.

Day's convolution
\[
(A \boxtimes B)(n,\beta) := \bigoplus_{\alpha + \gamma = \beta} \bigoplus_{i+j=n} \ind_{\sym_i \times \sym_j}^{\sym_n} (A(i,\alpha) \otimes B(j,\gamma))
\]
defines a symmetric monoidal structure on $[\sym \times \Upsilon,\ms{A}]$, with the identity object $S_0$. This induces a commutative ring structure on $\geKab$, whose product is defined on generators by
\[
[A] [B] := [A \boxtimes B] \,.
\]

For $\Upsilon = 0$, the Grothendieck ring $\eKab$ is isomorphic to the completion of $\Lambda \otimes K_0(\ms{A})$, where $\Lambda$ is the filtered ring of symmetric functions. This isomorphism sends $[s_n] \in \eKab$ to $s_{(n)} \otimes 1$, where $s_{(n)} = h_n$ is the homogeneous symmetric function of degree $n$, i.e., the Schur function corresponding to the partition $(n) \vdash n$. As a consequence, there is a ring isomorphism
\[
\eKab \cong \Kab[\![s_1,s_2,\dots]\!] \,.
\]

\begin{rmk}
In fact, $\eKab$ is isomorphic to the completion of $\Lambda \otimes K_0(\ms{A})$ as a $\lambda$-ring (see \cite[Thm.\ 4.8]{Get95conf}).
\end{rmk} 

Let $\gKab \cong \Kab[\![\Upsilon]\!]$ be the subring of $\geKab$ generated by the elements $[A]$ such that $A(n,\beta) = 0$ for all $n>0$. Its inclusion into $\geKab$ induces a commutative $\gKab$-algebra structure on $\geKab$. The algebra $\geKab$ is filtered by the subalgebras
\[
F_n := \left\langle \bigcup_{i=0}^n \bigl\{ [A] \mid A(j,\beta) = 0\ \text{if}\ j<n-i\ \text{or}\ m_\beta < i \bigr\} \right\rangle ,
\]
and it is complete with respect to this filtration. 

The composition (usually known as plethysm in this context)
\[
(A \circ B)(n,\beta) := \bigoplus_{\alpha + \gamma = \beta} \bigoplus_{i \ge 0} A(i,\alpha) \otimes_{\sym_i} B^{\boxtimes i}(n,\gamma) \,,
\]
which is defined whenever $B(0,0) = 0$, is the tensor product of a partially defined nonsymmetric monoidal structure on $[\sym \times \Upsilon, \ms{A}]$. The identity object of this monoidal structure is $S_1$. The functor $- \circ B$ preserves biproducts and commutes with Day's convolution, whereas this does not hold true for the functor $A \circ -$. For example, if $A = S_n$ then
\begin{equation} \label{schurab}
S_n \circ (B \oplus B') \cong \bigoplus_{i+j=n} (S_i \circ B) \boxtimes (S_j \circ B') \,.
\end{equation}

By means of \eqref{schurab} and the isomorphism $\geKab \cong (\Kab[\![s_1,s_2,\dots,]\!])[\![\Upsilon]\!]$, one easily sees that the composition induces a unique associative operation $\circ: \geKab \times F_1 \to \geKab$ such that
\[
[A] \circ [B] := [A \circ B] \,.
\]
This operation is a $\gKab$-algebra homomorphism in the first argument. With respect to the second argument, it is conveniently expressed via the elements $p_n \in \geKab$, which are analogous to the power sum symmetric functions in $\Lambda$. These elements are recursively defined by
\[
p_1 = s_1 \,, \qquad n s_n = p_n + s_1 p_{n-1} + \dots + s_{n-1} p_n \,,
\]
and there is an isomorphism
\[
\geKab \otimes_{\Z} \Q \cong \bigl( (\Kab \otimes_\Z \Q)[\![p_1,p_2,\dots]\!]) \bigr) [\![\Upsilon]\!] \,.
\]
The convenience of considering the elements $p_n$ lies in the next result, which follows from \cite[\S I.8]{Mac95}.

\begin{prop} \label{Newtonpoly}
For any $n \in \Z_{>0}$, the map $p_n \circ -$ is a $\gKab$-algebra homomorphism. Moreover, $p_n \circ p_m = p_{nm}$ for all $n,m$.
\end{prop}

As we will see in \S\ref{esempiacaso}, Proposition \ref{Newtonpoly} provides a way to effectively compute $a \circ b$ for any $a \in \geKab$ and $b \in F_1$.

Now, in order to define a composition algebra structure on $\geKab$, the last item we need is a suitable derivation. As in the case of $\ms{Var}$, we can consider the endofunctor $D$ of $[\sym \times \Upsilon, \ms{A}]$ given by
\[
(DA)(n,\beta) := A(n + 1, \beta) \,.
\]
Since this endofunctor preserves biproducts and satisfies both Leibniz's rule
\[
D(A \boxtimes B) \cong (DA \boxtimes B) \oplus (A \boxtimes DB)
\]
and the chain rule
\[
D(A \circ B) \cong (DA \circ B) \boxtimes DB \,,
\]
it induces a $\gKab$-derivation $D$ on $\geKab$, which is defined on generators by $D[A] := [DA]$. When applied to the elements $s_n$ and $p_n$, the derivation $D$ has a simple form: for all $n \ge 1$, $Ds_n = s_{n-1}$ and $Dp_n = \delta_{n,1}$, where $\delta_{i,j}$ is the Kronecker delta.

As a result of the preceding discussion, the following theorem is obtained (see \cite[Thm.\ 5.1]{GP06}).

\begin{thm}
The Grothendieck group $\geKab$ is a complete composition algebra over $\gKab$.
\end{thm}

The category of mixed Hodge structures over $\Q$ is a $\Q$-linear rigid abelian tensor category, thus all the above results can be applied to $\ms{A} = \ms{MHS}$. In view of the purpose of this paper, we will be specifically concerned with the case where $\Upsilon = \Z_{\ge 0}$, in which we have
\[
K_0(\ms{MHS}^{\Z_{\ge 0}}) \cong \Kmhs\ser \qquad \text{and} \qquad K_0^{\sym}(\ms{MHS}^{\Z_{\ge 0}}) \cong \eKmhs\ser \,.
\] 

\begin{rmk}
The Grothendieck group $\Kmhs$ is in fact isomorphic to the Grothendieck group $K_0(\ms{HS})$ of pure Hodge structures over $\Q$ (see for instance \cite[Cor.\ 3.9]{PS08}). The isomorphism $\Kmhs \to K_0(\ms{HS})$ is defined on generators by $[(V,W_\bullet,F^\bullet)] \mapsto \sum_{i \in \Z} [\mr{Gr}^W_i(V)]$. In particular, this shows that the information about the weight filtration is lost in $\Kmhs$.
\end{rmk}

\subsection{The Hodge-Grothendieck characteristic}

Let $X$ be a variety over $\C$. By Deligne's theory of mixed Hodge structures, each compactly supported cohomology group $H^i_\mr{c}(X,\Q)$ carries a natural mixed Hodge structure $(H^i_\mr{c}(X,\Q), W_{\bullet}, F^{\bullet})$ over $\Q$.

\begin{de} \label{hgcharsenzaz}
The Hodge-Grothendieck characteristic (hereinafter referred to as HG-characteristic) $\hg: \Kvar \to \Kmhs$ is the ring homomorphism that is defined on generators by
\[
\hg(X) := \sum_{i \ge 0} (-1)^i \bigl[\bigl( H^i_\mr{c}(X,\Q), W_{\bullet}, F^{\bullet} \bigr)\bigr] \,.
\]
Here and henceforth, $\hg(X)$ is a shortened notation for $\hg([X])$.
\end{de}

\begin{rmk}
The homomorphism $\hg$ has different names in the literature. We follow the terminology of \cite{PS08}. In \cite{GP06}, $\hg$ is called Serre characteristic. The proof of the existence of $\hg$ and its properties first appeared in \cite{DK86}.
\end{rmk}

For instance, the HG-characteristic of $\maff = [\aff^1_\C]$ is $[\Q(-1)] \in \Kmhs$. This class is denoted by $\caff$. 

The HG-characteristic is a refined version of the $E$-polynomial, which is indeed the composition of $\hg: \Kvar \to \Kmhs$ with the ring homomorphism $\Kmhs \to \Z[t,u,t^{-1},u^{-1}]$ given by
\[
[(V,W_{\bullet},F^{\bullet})] \mapsto \sum_{p,q \in \Z} \dim_\C\bigl( \mr{Gr}_F^p \mr{Gr}^W_{p+q}(V \otimes_\Q \C) \bigr) t^p u^q \,.
\]

If the variety $X$ carries an action of $\sym_n$, then the $\sym_n$-representation $H^i_\mr{c}(X,\Q)$ is compatible with the mixed Hodge structure on it, i.e., it is a $\sym_n$-representation in the category $\ms{MHS}$. This shows that Definition \ref{hgcharsenzaz} can be adapted to the context of $\Upsilon$-graded $\sym$-varieties. 

\begin{de}
The HG-characteristic $\hg: \geKvar \to \geKmhs$ is the ring homomorphism which is defined on generators by
\[
\hg(X) := \sum_{i \ge 0} (-1)^i \bigl[ H^i_\mr{c}(X,\Q) \bigr] \,,
\]
where $H^i_\mr{c}(X,\Q)(n,\beta) := (H^i_\mr{c}(X(n,\beta),\Q), W_{\bullet}, F^{\bullet})$ for all $n \in \Z_{\ge 0}$ and $\beta \in \Upsilon$. 
\end{de}

One of the main properties of the HG-characteristic is that it preserves the structures we have introduced in the previous sections. More precisely, we have the following result (see \cite[\S 5]{GP06}).

\begin{prop} \label{johnny}
The HG-characteristic $\hg: \geKvar \to \geKmhs$ is a homomorphism of complete composition algebras.
\end{prop}

For any projective variety $X$ that has at most finite quotient singularities, the knowledge of $\hg(X)$ determines its Hodge numbers, because the mixed Hodge structure $(H^i(X,\Q), W_{\bullet}, F^{\bullet})$ on its $i$-th cohomology group is in fact pure of weight $i$ (see \cite[Thm\ 2.43]{PS08}). In particular, this applies to the moduli spaces of genus $0$ stable maps to a finite dimensional Grassmannian. Therefore, the next sections will be concerned with the computation of the HG-characteristic of these spaces.


\section{Motivic study of spaces of genus $0$ stable maps} \label{echenneso}

Let $Y$ be a nonsingular, projective, convex variety, and let $\Upsilon \subseteq H_2(Y,\Z)$ be its monoid of curve classes. For any $n \in \Z_{\ge 0}$ and $\beta \in \Upsilon$, the coarse moduli space $\cls{M}_{0,n}(Y,\beta)$ of stable $(Y,0,n,\beta)$-maps is a normal projective variety with finite quotient singularities \cite[Thm.\ 2]{FP97}.

The variety $\cls{M}_{0,n}(Y,\beta)$ admits a stratification by locally closed subvarieties indexed by stable $(n,\beta)$-trees. When $Y$ is a finite-dimensional Grassmanian, this stratification plays a key role in the development of a recursive procedure to compute the HG-characteristic of $\cls{M}_{0,n}(Y,\beta)$.

In this section, we follow the same approach as that adopted in \cite[\S 4]{GP06} for the case $Y = \pr^r$.

\subsection{Stratification by stable marked trees}

In order to fix our notation, we first recall some definitions regarding graphs.

\begin{de}
A graph $\tau$ is a quadruple $(F_\tau, V_\tau, \partial_\tau, j_\tau)$, where $F_\tau$ is a finite set of flags, $V_\tau$ is a finite set of vertices, $\partial_\tau: F_\tau \to V_\tau$ is a map, and $j_\tau: F_\tau \to F_\tau$ is an involution. The elements of $E_\tau := \{ \{f,f'\} \subseteq F_\tau \mid j_\tau(f) \neq f'\}$ $F_\tau$ are called the edges of $\tau$, and the elements of $L_\tau := \{ f \in F_\tau \mid j_\tau(f) = f \}$ are called the leaves (or tails) of $\tau$. For each vertex $v \in V_\tau$, its valence $n(v)$ is the cardinality of $F_\tau(v) := \partial_\tau^{-1}(\{v\})$.

An isomorphism of graphs $\phi: \tau \iso \sigma$ is a pair $(\phi_F, \phi_V)$, where $\phi_F: F_\tau \to F_\sigma$ and $\phi_V: V_\tau \to V_\sigma$ are bijective maps such that the diagrams
\[
\begin{tikzcd}
F_\tau \arrow{r}{\phi_F} \arrow{d}{\partial_\tau} & F_\sigma \arrow{d}{\partial_\sigma} \\
V_\tau \arrow{r}{\phi_V} & V_\sigma
\end{tikzcd}
\qquad \text{and} \qquad
\begin{tikzcd}
F_\tau \arrow{r}{\phi_F} \arrow{d}{j_\tau} & F_\sigma \arrow{d}{j_\sigma} \\
F_\tau \arrow{r}{\phi_F} & F_\sigma
\end{tikzcd}
\]
commute.
\end{de}

\begin{de}
A tree is a graph $\tau$ such that 
\begin{enumerate}[label=(\roman*)]
\item $\#V_\tau = \#F_\tau + 1$, and 
\item for any $v,v' \in V_\tau$, there exists a sequence of edges connecting $v$ and $v'$, i.e., there exist $f_1,\dots,f_{2k} \in F_\tau$ such that $\partial_\tau(f_1)=v$, $\partial_\tau(f_{2k})=v'$, $j_\tau(f_{2i+1}) = f_{2i+2}$ for all $i=0,\dots,k-1$, and $\partial_\tau(f_{2i}) = \partial_\tau(f_{2i+1})$ for all $i=1,\dots,k-1$.
\end{enumerate}
In other words, a tree is a graph whose geometric realization is simply connected (see \cite[Def.\ 1.2]{BM96} for the geometric realization of a graph).
\end{de}

Now, let $\Upsilon$ be the monoid of curve classes of a nonsingular, projective, convex variety $Y$. More generally, one can replace $\Upsilon$ by any semigroup with indecomposable zero. We consider trees that are marked by elements of $\Upsilon$, and whose leaves are labelled by nonnegative integers.

\begin{de}
Let $n \in \Z_{\ge 0}$ and $\beta \in \Upsilon$. An $(n,\beta)$-tree is a triple $(\tau, l_\tau, \beta_\tau)$, where $\tau$ is a tree, $l_\tau: L_\tau \to \{1,\dots,n\}$ is a bijective map, and $\beta_\tau: V_\tau \to \Upsilon$ is a map such that $\sum_{v \in V_\tau} \beta_\tau(v) = \beta$.

An isomorphism of $(n,\beta)$-trees $(\tau, l_\tau, \beta_\tau) \iso (\sigma, l_\sigma, \beta_\sigma)$ is an isomorphism $\phi: \tau \iso \sigma$ of the underlying graphs such that the diagrams
\[
\begin{tikzcd}[column sep=small]
L_\tau \arrow{rr}{\phi_F|_{L_\tau}} \arrow[swap]{rd}{l_\tau} & & L_\sigma \\
& \{1,\dots,n\} \arrow[leftarrow, swap]{ru}{l_\sigma} &
\end{tikzcd}
\qquad \text{and} \qquad
\begin{tikzcd}
V_\tau \arrow{rr}{\phi_V} \arrow[swap]{rd}{\beta_\tau} & & V_\sigma \\
& \Upsilon \arrow[leftarrow, swap]{ru}{\beta_\sigma} &
\end{tikzcd}
\]
commute.
\end{de}

For each stable $(Y,0,n,\beta)$-map $(C,x,f)$, there is an associated $(n,\beta)$-tree $(\tau, l_\tau, \beta_\tau)$, which is defined as follows.
\begin{itemize}
\item Let $\nu: C^\nu \to C$ denote the normalization of $C$. The flags of $\tau$ correspond to the closed points $y \in C^\nu$ such that $\nu(y)$ is a special point of $(C,x)$.
\item The vertices of $\tau$ correspond to the irreducible components of $C$, and $\partial_\tau: F_\tau \to V_\tau$ sends $y$ to the component where $\nu(y)$ lies.
\item If $\nu(y)$ is a marked point, then $j_\tau(y) := y$. If $\nu(y)$ is a double point, then $j_\tau(y) := y'$, where $y' \neq y$ is the other closed point of $C^\nu$ such that $\nu(y) = \nu(y')$. In particular, the edges of $\tau$ correspond to the nodes of $C$, and the leaves of $\tau$ correspond to the marked points of $C$.
\item The labelling $l_\tau: L_\tau \to \{1,\dots,n\}$ of the leaves is determined by the markings $x = (x_i)_{1 \le i \le n}$: $l_\tau(y) := i$ if and only if $\nu(y) = x_i$.
\item The map $\beta_\tau: V_\tau \to \Upsilon$ sends each irreducible component $C'$ of $C$ to the class of $f|_{C'}$.
\end{itemize}
The stability condition on $(C,x,f)$ corresponds to a stability condition on its dual $(n,\beta)$-tree.

\begin{de}
An $(n,\beta)$-tree $\tau$ is stable if for each $v \in V_\tau$ either $\beta_\tau(v) \neq 0$ or $n(v)>2$. 
\end{de}

In particular, a $(Y,0,n,\beta)$-map is stable if and only if its dual $(n,\beta)$-tree is stable. Furthermore, note that the dual $(n,\beta)$-trees of two isomorphic stable $(Y,0,n,\beta)$-maps are isomorphic as well.

Conversely, given an isomorphism class $[\tau]$ of stable $(n,\beta)$-trees, one can consider the locus $M(\tau)$ in $\cls{M}_{0,n}(Y,\beta)$ that parametrizes those maps whose dual $(n,\beta)$-tree is isomorphic to $\tau$. This locus is a locally closed subvariety of $\cls{M}_{0,n}(Y,\beta)$ of codimension $\#E_\tau$ (see \cite[\S 6]{FP97}). 

For instance, the locus that corresponds to the isomorphism class of stable $(n,\beta)$-trees with a single vertex and $n$ leaves is the open subvariety $M_{0,n}(Y,\beta)$, which parametrizes maps from nonsingular curves. If $[\tau]$ is a different isomorphism class, then $M(\tau)$ lies in the boundary $\cls{M}_{0,n}(Y,\beta) \setminus M_{0,n}(Y,\beta)$.

The decomposition of $\cls{M}_{0,n}(Y,\beta)$ into the locally closed subvarieties $M(\tau)$ is finite, as the next result shows.

\begin{lmm} \label{alberisonofiniti}
The set $\Gamma_{0,n}(\beta)$ of isomorphism classes of stable $(n,\beta)$-trees is finite.
\end{lmm}

\proof
The proof is the same as that of \cite[Prop.\ 4.3]{GP06} for $Y=\pr^r$, adapted to our more general case. For every $(n,\beta)$-tree $\tau$, we have
\[
\#V_\tau = \#E_\tau + 1 = \frac{1}{2} \Biggl( \sum_{v \in V_\tau} n(v) - n \Biggr) + 1 \,,
\]
thus $\sum_{v \in V_\tau} (n(v)-2) = n-2$. In particular, the number of vertices $v \in V_\tau$ such that $n(v)>2$ is bounded by $n-2$. Now, the set 
\[
\left\{ (\beta_1,\dots,\beta_i \in \Upsilon^i \mid \beta_1+\dots+\beta_i=\beta \,,\ \beta_j \neq 0\ \forall\,j \right\}
\]
is empty for almost all $i$; let $i_\beta$ be the greatest integer for which it is nonempty. If $\tau$ is stable, then the number of vertices $v \in V_\tau$ such that $n(v) \le 2$ is bounded by $i_\beta$. Therefore, we have
\[
\#F_\tau = \sum_{v \in V_\tau} n(v) = n-2 + 2(\#V_\tau) \le 3n + i_\beta - 4 \,.
\]
Since the number of isomorphism classes of trees with a fixed number of flags is finite, it follows that $\Gamma_{0,n}(\beta)$ is a finite set.
\endproof

As a consequence of Lemma \ref{alberisonofiniti}, the equality
\[
\bigl[\cls{M}_{0,n}(Y,\beta)\bigr] = \sum_{[\tau] \in \Gamma_{0,n}(\beta)} [ M(\tau) ]
\]
holds in $\Kvar$. 

The variety $\cls{M}_{0,n}(Y,\beta)$ has a natural left action of the symmetric group $\sym_n$, which is given by permutation of the marked points. If we consider its class in $\nKvar{n}$, then in general the above equation does not hold, because the loci $M(\tau)$ may be not $\sym_n$-invariant. Only the weaker equalies
\[
\bigl[\cls{M}_{0,n}(Y,\beta)\bigr] = \left[ \coprod_{[\tau] \in \Gamma_{0,n}(\beta)} M(\tau) \right] = [M_{0,n}(Y,\beta)] + \left[ \coprod_{[\tau] \in \Gamma_{0,n}(\beta),\, E_\tau \neq \emptyset} M(\tau) \right]
\]
hold in $\nKvar{n}$.

As a consequence of \cite[\S 6]{FP97} (see also \cite[\S XII.10]{ACG11} for the analogous description in the case of moduli of pointed curves), the loci $M(\tau)$ can be given a description in terms of the gluing of stable maps from nonsingular curves. For this purpose, it is useful to work with stable maps from pointed curves whose marked points are labelled by an arbitrary finite set.

\begin{de}
For any finite set $L$ and any $\beta \in \Upsilon$ such that either $\beta \neq 0$ or $\#L > 2$, $M_{0,L}(Y,\beta)$ is the coarse moduli space that parametrizes equivalence classes of triples $(C,\iota,f)$, where $C$ is a nonsingular projective curve of genus $0$, $\iota: L \to C$ is a closed immersion, and $f: C \to Y$ is a morphism such that $f_\ast([C]) = \beta$.
\end{de} 

The variety $M_{0,L}(Y,\beta)$ is (non-canonically) isomorphic to $M_{0,\#L}(Y,\beta)$, and it is equipped with a canonical evaluation morphism $M_{0,L}(Y,\beta) \to Y^L$, which maps each geometric point $[(C,\iota,f)]$ to $f \circ \iota: L \to Y$. 

Let $\tau$ be a stable $(n,\beta)$-tree. Since $F_\tau = \bigsqcup_{v \in V_\tau} F_\tau(v)$, the evaluation morphisms $M_{0,F_\tau(v)}(Y,\beta_\tau(v)) \to Y^{F_\tau(v)}$ induce a morphism
\[
\prod_{v \in V_\tau} M_{0, F_\tau(v)}(Y,\beta_\tau(v)) \to Y^{F_\tau} \,.
\]
There is also a natural injection $Y^{E_\tau \sqcup L_\tau} \to Y^{F_\tau}$, which is given by composing with the surjective map $F_\tau \to E_\tau \sqcup L_\tau$ that maps each flag to its orbit under $j_\tau$. Let $M_\Box(\tau)$ be the fiber product
\[
\begin{tikzcd}
M_\Box(\tau) \arrow{r} \arrow{d} & \prod_{v \in V_\tau} M_{0, F_\tau(v)}(Y,\beta_\tau(v)) \arrow{d} \\
Y^{E_\tau \sqcup L_\tau} \arrow{r} & Y^{F_\tau}
\end{tikzcd} 
\]
with respect to these morphisms. The next result follows from \cite[\S 6.2]{FP97} (compare also with \cite[Prop.\ 10.11]{ACG11} in the case of moduli of curves).

\begin{prop} \label{taufibprod}
For each $[\tau] \in \Gamma_{0,n}(Y,\beta)$, there is a canonical isomorphism 
\[
M(\tau) \cong M_\Box(\tau)/\Aut(\tau) \,.
\]
\end{prop}

\proof
By \cite[Prop.\ 2.4]{BM96} and the definition of $M_\Box(\tau)$, there is a canonical gluing morphism $M_\Box(\tau) \to M(\tau)$. This morphism is invariant under automorphisms of $\tau$, therefore it factors through $M_\Box(\tau)/\Aut(\tau)$. The resulting morphism $M_\Box(\tau)/\Aut(\tau) \to M(\tau)$ is the normalization of $M(\tau)$, which is an isomorphism because $M(\tau)$ is already normal.
\endproof

\subsection{Stable maps from genus $0$ curves to Grassmannians} \label{sectonrecrel}

Let $V$ be a fixed $\C$-vector space of dimension $k \in \Z_{> 0}$, and let $r$ be an integer such that $0 < r < k$. We consider stable maps with values in the Grassmannian $\grass := G(r,V)$ of $r$-dimensional quotients of $V$. 

Under the isomorphism $H_2(\grass,\Z) \cong \Z$, the monoid of curve classes of $\grass$ is identified with $\Z_{\ge 0}$. Therefore, a curve class in $\grass$ will be denoted by $d \in \Z_{\ge 0}$.

\begin{prop} \label{evalloctriv}
The evaluation morphisms $\mr{ev}_a: \cls{M}_{g,n}(\mathds{G}, d) \to \mathds{G}$ are locally trivial fibrations in the Zariski topology.
\end{prop}

\proof
For notational simplicity, we identify $V$ with $\C^k$, via the choice of a basis of $V$. Throughout all the proof, by a point we mean a closed point. For any $k \times k$ matrix $A$ and any $I,J \subseteq \{1,\dots,k\}$, we denote by $A_{(I,J)}$ the submatrix of $A$ whose rows are indexed by $I$ and whose columns are indexed by $J$.

For every $I= \{i_1,\dots,i_{k-r}\} \subseteq \{1,\dots,k\}$, let $U_I \cong \aff^{r(k-r)}$ be the open affine subset of $\grass$ whose points are the quotients $[q: \C^k \to W]$ such that $\ker(q) \cap \langle \{e_j \mid j \notin I\} \rangle = \{0\}$. If $I^c=\{1,\dots,k\} \setminus I$, then $U_I$ is isomorphic to the closed subgroup $H_I$ of $\GL(k,\C)$ consisting of those matrices $A \in \GL(k,\C)$ such that $A_{(I,I)} = \1_{k-r}$, $A_{(I^c,I^c)} = \1_r$ and $A_{(I,I^c)}= 0$. An isomorphism $\phi: U_I \to H_I$ is defined as follows: for any $[q: \C^k \to W] \in U_I$, let $K_q$ be the $k \times (k-r)$ matrix whose column vectors span $\ker(q)$ and whose $I$-th minor is $\1_{k-r}$; then $\phi(q) \in H_I$ is the matrix such that $\phi(q)_{(I^c,I)}$ is the $I^c$-th minor of $K_q$. The restriction of the $\GL(k,C)$-action on $\grass$ yields an action of $H_I$ on $\grass$ such that $U_I$ is $H_I$-invariant. The $H_I$-action on $U_I$ corresponds to the translation in $\aff^{r(k-r)}$; in particular, it is free and transitive. 

Let $F$ denote the fiber $(\mr{ev}_a)^{-1}(0_I)$, where $0_I \in U_I$ is the point corresponding to $0 \in \aff^{r(k-r)}$. The isomorphism class of $F$ does not depend on the choice of $I$: for any other $J \subseteq \{1,\dots,k\}$ with $\#J=k-r$, the action of row-switching matrices of $\GL(k,\C)$ on $\grass$ determines an isomorphism $F \cong (\mr{ev}_a)^{-1}(0_J)$. Since the $U_I$'s form an open cover of $\grass$, the proposition will be proved if we show that there is an isomorphism $(\mr{ev}_a)^{-1}(U_I) \cong U_I \times F$.

Let $(\pi: C \to T,x,f)$ be a stable $(\grass, g, n, d)$-map over $T$ such that $f \circ x_a$ factors through $U_I$. Then we can consider the morphism $f_0: C \to \grass$ that sends $z \in C$ to $\phi((f \circ x_a \circ \pi)(z))^{-1} \cdot f(z)$; informally, we are translating $f \circ x_a$ to $0_I$. The map $f_0$ is a morphism of class $d$ such that $f_0 \circ x_a$ factors through $\{0_I\}$. As a consequence, the mapping
\[
(\pi: C \to T,x,f) \mapsto (f \circ x_a: T \to U_I, (\pi: C \to T, x, f_0))
\]
determines a canonical morphism $\psi: (\mr{ev}_a)^{-1}(U_I) \to U_I \times F$.

Conversely, let $(\pi': C' \to T, x', f')$ be a stable $(\grass,g,n,d)$-map over $T$ such that $f' \circ x'_l$ factors through $\{0_I\}$, and let $h: T \to U_I$ be a morphism. Then the map $f_h: C \to \grass$ defined by $f_h(z) := \phi((h \circ \pi')(z)) \cdot f'(z)$ is a morphism of class $d$ such that $f_h \circ x'_a = h$. Therefore, the mapping
\[
(h: T \to U_I, (\pi': C' \to T, x', f')) \mapsto (\pi': C' \to T, x', f_h)
\]
determines a canonical morphism $\psi': U_I \times F \to (\mr{ev}_a)^{-1}(U_I)$, which is clearly the inverse of $\psi$.
\endproof

The fibers of $\mr{ev}_{n+1}: \cls{M}_{0,n+1}(\grass, d) \to \grass$ and of its restriction to $M_{0,n+1}(\grass,d)$ will play a special role in our treatment.

\begin{de}
If $p \in \grass$ is any closed point, $\cls{\Phi}_{n,d}$ is the fiber of the morphism $\mr{ev}_{n+1}: \cls{M}_{0,n+1}(\grass, d) \to \grass$ over $p$. Similarly, $\Phi_{n,d}$ is the fiber of $\mr{ev}_{n+1}|_{M_{0,n+1}(\grass,d)}$ over $p$.
\end{de}

Since we are only interested in the isomorphism class of the fibers, the choice of $p \in \grass$ in the above definition is arbitrary.

By identifying $\sym_n$ with $\{ \sigma \in \sym_{n+1} \mid \sigma(n+1)=n+1 \}$, the restriction of the $\sym_{n+1}$-action on $\cls{M}_{0,n+1}(\grass, d)$ (resp.\ $M_{0,n+1}(\grass, d)$) yields an $\sym_n$-action on $\cls{\Phi}_{n,d}$ (resp.\ $\Phi_{n,d}$). Let $\cls{M}$, $M$, $\cls{\Phi}$ and $\Phi$ be the following graded $\sym$-varieties:
\begin{align*}
\cls{M}(n,d) &:= \cls{M}_{0,n}(\grass), d) \,, & M(n,d) &:= M_{0,n}(\grass, d) \,, \\
\cls{\Phi}(n,d) &:= \cls{\Phi}_{n,d} \,, & \Phi(n,d) &:= \Phi_{n,d} \,.
\end{align*}
The next result is a direct consequence of Proposition \ref{evalloctriv}.

\begin{cor} \label{derivserre}
The equalities
\[
D[\cls{M}] = [\grass] [\cls{\Phi}] \qquad \text{and} \qquad D[M] = [\grass] [\Phi]
\]
hold in $\eKvar\ser$. 
\end{cor}

Following a procedure analogous to \cite[\S 4]{GP06}, we are going to find a formula that relates the classes of $\cls{M}$, $M$, $\cls{\Phi}$ and $\Phi$ in $\eKvar\ser$. First, we need some preliminary results.

Let us fix a representative for each isomorphism class in $\Gamma_{0,n}(d)$. The set of these representatives is denoted by $\tilde{\Gamma}_{0,n}(d)$. For any $\tau \in \tilde{\Gamma}_{0,n}(d)$, let us also fix representatives for each equivalence class in $V_\tau/\Aut(\tau)$, $E_\tau/\Aut(\tau)$ and $F_\tau/\Aut(\tau)$; the corresponding sets are denoted by $\tilde{V}_\tau$, $\tilde{E}_\tau$ and $\tilde{F}_\tau$, respectively. Notice that $L_\tau/\Aut(\tau) = L_\tau$. 
For any $a \in V_\tau \sqcup E_\tau \sqcup L_\tau \sqcup F_\tau$, $\Aut(\tau,a)$ denotes the group of automorphisms of $\tau$ that fix $a$.

Finally, let us introduce the sets
\begin{align*}
\Gamma_{0,n}^V(d) &:= \{ (\tau,v) \mid \tau \in \tilde{\Gamma}_{0,n}(d)\,, v \in \tilde{V}_\tau \} \,, \\
\Gamma_{0,n}^E(d) &:= \{ (\tau,e) \mid \tau \in \tilde{\Gamma}_{0,n}(d)\,, e \in \tilde{E}_\tau \} \,, \\
\Gamma_{0,n}^L(d) &:= \{ (\tau,l) \mid \tau \in \tilde{\Gamma}_{0,n}(d)\,, l \in L_\tau \} \,, \\
\Gamma_{0,n}^F(d) &:= \{ (\tau,\xi) \mid \tau \in \tilde{\Gamma}_{0,n}(d)\,, \xi \in \tilde{F}_\tau \} \,.
\end{align*}
The following result translates the combinatorial properties of the moduli spaces we are considering into equations in $\eKvar\ser$ (cf.\ \cite[Lemma 4.7]{GP06}).

\begin{prop} \label{fundlemma}
The equations
\begin{align}
\sum_{d=0}^\infty q^d \sum_{n=0}^\infty \left[ \coprod_{(\tau,v) \in \Gamma_{0,n}^V(d)} M_\Box(\tau)/\Aut(\tau,v) \right] &= [M] \circ (s_1 + [\cls{\Phi}] ) \,, \label{vertex} \\
\sum_{d=0}^\infty q^d \sum_{n=0}^\infty \left[ \coprod_{(\tau,e) \in \Gamma_{0,n}^E(d)} M_\Box(\tau)/\Aut(\tau,e) \right] &= [\grass] (s_2 \circ [\cls{\Phi}]) \,, \label{edge} \\
\sum_{d=0}^\infty q^d \sum_{n=0}^\infty \left[ \coprod_{(\tau,l) \in \Gamma_{0,n}^L(d)} M_\Box(\tau)/\Aut(\tau,l) \right] &= s_1 D[\cls{M}] \,, \label{leaf} \\
\sum_{d=0}^\infty q^d \sum_{n=0}^\infty \left[ \coprod_{(\tau,\xi) \in \Gamma_{0,n}^F(d)} M_\Box(\tau)/\Aut(\tau,\xi) \right] &= [\grass] [\cls{\Phi}]^2 + s_1 D[\cls{M}] \label{flag}
\end{align}
hold in $\eKvar\ser$.
\end{prop}

\proof
The result is a direct consequence of Proposition \ref{taufibprod}, Proposition \ref{evalloctriv} and the definition of $\circ$. For the sake of geometric intuition, a geometric interpretation of the equalities \eqref{vertex}--\eqref{flag} is provided below. In what follows, $(C,x,f)$ denotes a stable $(\grass,0,n,d)$-map with dual $(n,d)$-tree $\tau$. 

\eqref{vertex} If $v \in V_\tau$ is a fixed vertex and $C_v$ is the corresponding irreducible component of $C$, let $\nu: \tilde{C}_v \to C_v$ be the normalization of $C_v$, let $z_i$, for $i \in F_\tau(v)$, be the closed points of $\tilde{C}_v$ that lie over the special points of $C_v$, and let $f_v = f|_{C_v} \circ \nu$. Then $(C,x,f)$ is obtained by gluing the stable map $(\tilde{C}_v, (z_i)_{i \in F_\tau(v)}, f_v)$ with
\begin{itemize}
\item a stable map at each $z_i$ such that $\nu(z_i)$ is a node, and
\item a point at each $z_i$ such that $\nu(z_i)$ is a marked point.
\end{itemize}

\eqref{edge} If $e \in E_\tau$ is a fixed edge and $x_e$ is the corresponding node of $C$, then $(C,x,f)$ is obtained by gluing two stable maps, where both of them have an additional marked point that is mapped to $f(x_e) \in \grass$, at these additional points.

\eqref{leaf} Let $l \in L_\tau$ be a fixed leaf, and let $\lambda \in \{1,\dots,n\}$ be its label. Then the action of $\sym_n$ that permutes $x_1,\dots,x_n$ is the same as that induced by the action of its subgroup $\{ \sigma \in \sym_n \mid \sigma(\lambda) = \lambda\} \cong \sym_{n-1}$ that permutes $x_1,\dots,\hat{x}_\lambda,\dots,x_n$.

\eqref{flag} Let $\xi \in F_\tau$ be a fixed flag, and let $z_\xi$ be the corresponding point of $\tilde{C}$, where $\nu: \tilde{C} \to C$ is the normalization. If $\nu(z_\xi)$ is a marked point of $(C,x)$, then we can argue as for \eqref{leaf}. This argument yields the second summand on the right-hand side of \eqref{flag}. 

If $\nu(z_\xi)$ is a node of $C$, then we can argue as for \eqref{edge}. However, while in the case of edges there may be automorphisms of $(\tau,e=\{\xi, j_\tau(\xi)\})$ that exchange $\xi$ and $j_\tau(\xi)$, such automorphisms are not in $\Aut(\tau,\xi)$. This explains the lack of the quotient by $\sym_2$ in the first summand on the right-hand side of \eqref{flag}, as compared with \eqref{edge}.
\endproof

In order to prove the above-mentioned formula, we need two other lemmas. The first one is \cite[Lemma\ 4.8]{GP06}.

\begin{lmm} \label{injgraph}
For every stable $(n,d)$-tree $\tau$, there is a canonical injective map $\iota: F_\tau \into V_\tau \sqcup E_\tau \sqcup L_\tau$. 
\end{lmm}

The proof of the second one is straightforward.

\begin{lmm} \label{inutilemma}
Let $X$ be a quasi-projective variety, and let $A$ be a finite set. Let $G$ be a finite group acting on both $X$ and $A$. Then any choice of representatives for the equivalence classes in $A/G$ determines an isomorphism
\[
(A \times X)/G \cong \coprod_{[a] \in A/G} X/G_a \,,
\]
where $G_a = \{g \in G \mid ga=a\}$ is the stabilizer subgroup of $G$ with respect to the chosen representative $a$ of $[a]$.
\end{lmm}

We can now prove the main theorem of this section (cf.\ \cite[Thm.\ 4.5]{GP06}).

\begin{thm} \label{recrelgrot}
The equality
\[
[ \cls{M} ] = [M] \circ ( s_1 + [ \cls{\Phi} ] ) + [G(r,V)] \bigl( s_2 \circ [ \cls{\Phi} ] - [ \cls{\Phi} ]^2 \bigr)
\]
holds in $\eKvar\ser$.
\end{thm}

\proof
By Lemma \ref{inutilemma}, we have the equalities
\begin{align*}
\left[ \coprod_{(\tau,v) \in \Gamma_{0,n}^V(d)} M_\Box(\tau)/\Aut(\tau,v) \right] &= \left[ \coprod_{[\tau] \in \Gamma_{0,n}(d)} (V_\tau \times M_\Box(\tau))/\Aut(\tau) \right] , \\
\left[ \coprod_{(\tau,e) \in \Gamma_{0,n}^E(d)} M_\Box(\tau)/\Aut(\tau,e) \right] &= \left[ \coprod_{[\tau] \in \Gamma_{0,n}(d)} (E_\tau \times M_\Box(\tau))/\Aut(\tau) \right] , \\
\left[ \coprod_{(\tau,l) \in \Gamma_{0,n}^L(d)} M_\Box(\tau)/\Aut(\tau,l) \right] &= \left[ \coprod_{[\tau] \in \Gamma_{0,n}(d)} (L_\tau \times M_\Box(\tau))/\Aut(\tau) \right] .
\end{align*}
in $\nKvar{n}$. 

The sum of the right-hand sides is equal to
\begin{equation} \label{labelacaso}
\left[ \coprod_{[\tau] \in \Gamma_{0,n}(d)} ((V_\tau \sqcup E_\tau \sqcup L_\tau) \times M_\Box(\tau)) / \Aut(\tau) \right] .
\end{equation}
The class \eqref{labelacaso} can be rewritten by means of Lemma \ref{injgraph}. Indeed, let us consider the canonical injection $\iota: F_\tau \into V_\tau \sqcup E_\tau \sqcup L_\tau$. Since $\tau$ is a tree, we have $\#F_\tau = \#(V_\tau \sqcup E_\tau \sqcup L_\tau)-1$, thus the injection $\iota$ induces a canonical bijection $F_\tau \sqcup \{\bigcdot\} \iso V_\tau \sqcup E_\tau \sqcup L_\tau$. As a consequence, \eqref{labelacaso} equals
\[
\left[ \coprod_{[\tau] \in \Gamma_{0,n}(d)} ((F_\tau \sqcup \{\bigcdot\}) \times M_\Box(\tau)) / \Aut(\tau) \right] ,
\]
which in turn equals
\[
\left[ \coprod_{(\tau,i) \in \Gamma_{0,n}^F(d)} M_\Box(\tau)/\Aut(\tau,i) \right] + \left[ \coprod_{[\tau] \in \Gamma_{0,n}(d)} M_\Box(\tau)/\Aut(\tau) \right] ,
\]
again by Lemma \ref{inutilemma}. Furthermore, Proposition \ref{taufibprod} implies that
\[
\left[ \coprod_{[\tau] \in \Gamma_{0,n}(d)} M_\Box(\tau)/\Aut(\tau) \right] = \left[ \coprod_{[\tau] \in \Gamma_{0,n}(d)} M(\tau) \right] = \bigl[ \cls{M}_{0,n}(\grass, d) \bigr] \,.
\]

Finally, the above sequence of equalities, together with Proposition \ref{fundlemma}, yields the equation
\[
[M] \circ (s_1 + [\cls{\Phi}]) + [\grass] (s_2 \circ [\cls{\Phi}]) + s_1 D[\cls{M}]
=
[\grass] [\cls{\Phi}]^2 + s_1 D[\cls{M}] + [ \cls{M} ]
\]
in $\eKvar\ser$, from which the theorem follows.
\endproof

\begin{cor} \label{calcfiber}
The equality
\[
[\grass] \left( [\cls{\Phi}] - [\Phi] \circ (s_1 + [\cls{\Phi}]) \right) = 0
\]
holds in $\eKvar\ser$.
\end{cor}

\proof
We apply the derivation $D$ to both sides of Theorem \ref{recrelgrot}. By Corollary \ref{derivserre}, the left-hand side becomes $D[\cls{M}] = [\grass] [\cls{\Phi}]$.

Since $D$ is additive, it can be applied separately to each summand on the right-hand side. Using the properties of $D$ and Corollary \ref{derivserre}, we get
\[
\begin{split}
D\bigl( [M] \circ (s_1 + [\cls{\Phi}]) \bigr) &= \left( D[M] \circ (s_1 + [\cls{\Phi}]) \right) D(s_1 + [\cls{\Phi}]) \\
&= \left( ([\grass] [\Phi]) \circ (s_1 + [\cls{\Phi}]) \right) \left( 1 + D[\cls{\Phi}] \right) \\
&= [\grass] \left( [\Phi] \circ (s_1 + [\cls{\Phi}]) \right) \left( 1 + D[\cls{\Phi}] \right) ,
\end{split}
\]
then
\[
\begin{split}
D\bigl( [\grass] (s_2 \circ [\cls{\Phi}]) \bigr) &= [\grass]\, D( s_2 \circ [\cls{\Phi}] ) \\
&= [\grass] \left( s_1 \circ [\cls{\Phi}] \right) D[\cls{\Phi}] \\
&= [\grass] [\cls{\Phi}]\, D[\cls{\Phi}] \,,
\end{split}
\]
and finally
\[
D\bigl( [\grass] [\cls{\Phi}]^2 \bigr) = 2 [\grass] [\cls{\Phi}]\, D[\cls{\Phi}] \,.
\]

The combination of all the above equalities yields
\[
[\grass] \left( 1 + D[\cls{\Phi}] \right) \left( [\cls{\Phi}] - [\Phi] \circ (s_1 + [\cls{\Phi}]) \right) = 0 \,.
\]
Since $1 + D[\cls{\Phi}]$ is not a zero divisor in $\eKvar\ser$, this implies the corollary.
\endproof

As the HG-characteristic $\hg: \eKvar\ser \to \eKmhs\ser$ is a morphism of complete composition algebras, Theorem \ref{recrelgrot} and Corollary \ref{calcfiber} also imply the following result.

\begin{cor} \label{veryimpcoroll}
The equalities
\[
\hg(\cls{M}) = \hg(M) \circ \bigl( s_1 + \hg(\cls{\Phi}) \bigr) + \hg(\grass) \bigl( s_2 \circ \hg(\cls{\Phi}) - \hg(\cls{\Phi})^2 \bigr)
\]
and 
\[
\hg(\cls{\Phi}) = \hg(\Phi) \circ \bigl( s_1 + \hg(\cls{\Phi}) \bigr)
\]
hold in $\eKmhs\ser$.
\end{cor}

The main upshot of Theorem \ref{recrelgrot} and its corollaries is that they provide a recursive procedure to compute the HG-characteristic of $\cls{M}_{0,n}(\grass,d)$ for all $n$ and $d$. 

The first equation of Corollary \ref{veryimpcoroll} shows that $\hg(\cls{M}_{0,n}(\grass,d)) \in \nKmhs{n}$ can be obtained from the following elements:
\begin{enumerate}[label=(\roman*)]
\item $\hg(\grass) \in \Kmhs$;
\item $\hg(\cls{\Phi}_{m,\delta}) \in \nKmhs{m}$ for all $\delta \le d$ and $m \le \begin{cases}
n & \text{if}\ \delta<d \\
n-2 & \text{if}\ \delta=d
\end{cases}$\,; \label{fibre}
\item $\hg(M_{0,m}(\grass,\delta)) \in \nKmhs{m}$ for all $\delta \le d$ and $m \le n + d - \delta$. \label{aperto}
\end{enumerate}

Given $V$ and $r$, $[\grass] \in \Kvar$ (and thus $\hg(\grass)$) is easily computable by means of the recursive equality
\[
[\grass] = [G(r,k)] = [G(r,k-1)] + \maff^{k-r} [G(r-1,k-1)] \,,
\]
starting from the base cases
\[
[G(1,k)] = [G(k-1,k)] = [ \pr^{k-1} ] = \sum_{i=0}^{k-1} \maff^i \,.
\] 

By applying Corollary \ref{derivserre}, the second equality of Corollary \ref{veryimpcoroll} becomes
\begin{equation} \label{recfibre}
\hg(\cls{\Phi}) = \frac{ D(\hg(M)) }{ \hg(\grass) } \circ \bigl( s_1 + \hg(\cls{\Phi}) \bigr) \,.
\end{equation}
Note that \eqref{recfibre} yields a recursive algorithm for computing $\hg(\cls{\Phi}_{m,\delta})$. In particular, it shows that $\hg(\cls{\Phi}_{m,\delta})$ is determined by $\hg(M_{0,i}(\grass,\epsilon))$ for $\epsilon \le \delta$ and $i \le m + \delta - \epsilon + 1$. Therefore, the elements \ref{fibre} above can be obtained from the elements \ref{aperto}.

In conclusion, we see that the computation of $\hg(\cls{M}_{0,n}(\grass,d))$ recursively boils down to the calculation of $\hg(M_{0,m}(\grass,\delta))$.

\subsection{The HG-characteristic of the open stratum} \label{carattserreaperto}

Motivated by the previous results, we shall now face the problem of computing the HG-characteristic of the open locus $M_{0,n}(\grass,d) \subseteq \cls{M}_{0,n}(\grass,d)$, which parametrizes those maps whose domain curve is isomorphic to $\pr^1$.

Let $\Mor_d(\pr^1, \grass)$ be the variety that parametrizes morphisms from $\pr^1$ to $\grass$ of class $d$, and let 
\[
\mr{F}(\pr^1,n) := \bigl\{ (x_1,\dots,x_n) \in \bigl(\pr^1\bigr)^n \mid x_i \neq x_j\ \forall\, i \neq j \bigr\} 
\] 
be the configuration space of $n$ distinct points in $\pr^1$. The symmetric group $\sym_n$ acts on $\mr{F}(\pr^1,n)$ by permuting the $n$-tuples of points. If $n+3d \ge 3$, then there is a canonical $\sym_n$-equivariant isomorphism
\begin{equation} \label{isoconf}
M_{0,n}(\grass,d) \cong \bigl( \Mor_d(\pr^1, \grass) \times \mr{F}(\pr^1,n) \bigr)/\Aut(\pr^1) \,,
\end{equation}
where $\Aut(\pr^1) = \mr{PGL}(2,\C)$ acts on both $\Mor_d(\pr^1, \grass)$ and $\mr{F}(\pr^1,n)$ via its action on $\pr^1$. In the remaining cases, namely for $d=0$ and $n \le 2$, the stability condition implies that $M_{0,n}(\grass,0) = \emptyset$.

Let $\mr{F}(\pr^1)$ be the $\sym$-variety defined as $\mr{F}(\pr^1)(n) := \mr{F}(\pr^1,n)$. The class of $\mr{F}(\pr^1)$ in $\eKvar$ can be expressed via the Exp and Log functions of \cite[\S 3]{GP06}.

\begin{prop}[{\cite[Thm.\ 3.2]{GP06}}]
The equality
\[
[ \mr{F}(\pr^1) ] = 1 + \Exp\bigl( [\pr^1] \Log(s_1) \bigr)
\]
holds in $\eKvar$.
\end{prop}

The corresponding formula for $\hg(\mr{F}(\pr^1)) \in \eKmhs$ can be simplified using the properties of $\eKmhs$.

\begin{cor}[{\cite[\S 5]{GP06}}] \label{configserrecar}
Let $\mu$ denote the M{\"o}bius function. Then the equality
\[
\hg(\mr{F}(\pr^1)) = (1+p_1) \prod_{n=1}^\infty (1+p_n)^{(1/n) \sum_{k|n} \mu(n/k) \caff^k}
\]
holds in $\eKmhs$. 
\end{cor}

The class of $\Aut(\pr^1)$ in $\Kvar$ and its HG-characteristic are also easily computable. Indeed, we have $[\mr{GL}(2,\C)] = \maff (\maff-1)^2 + \maff^2 (\maff-1)^2$, 
therefore
\[
[\Aut(\pr^1)] = [\mr{PGL}(2,\C)] = \frac{ [\mr{GL}(2,\C)] }{ [\C^\ast] } = \maff^3 - \maff
\]
and
\[
\hg(\Aut(\pr^1)) = \caff^3 - \caff \,.
\]

However, since the quotient morphism $\Mor_d(\pr^1, \grass) \times \mr{F}(\pr^1,n) \to M_{0,n}(\grass, d)$ is not a locally trivial fibration in the Zariski topology, Proposition \ref{solitarobatriviale} cannot be used to express $[M_{0,n}(\grass,d)] \in \nKvar{n}$ in terms of $[\mr{F}(\pr^1,n)]$ and $[\Aut(\pr^1)]$. This problem is overcome by passing to $\nKmhs{n}$ via the HG-characteristic.

\begin{thm}[{\cite[Thm.\ 5.4]{GP06}}] \label{GPtrivial}
Let $X$ be an $\sym$-variety, and let $G$ be a connected algebraic group which acts on each $X(n)$. Assume that the actions of $\sym_n$ and $G$ on each $X(n)$ commute. If the action of $G$ has finite stabilizers, then $\hg(X) = \hg(G)\, \hg(X/G)$ in $\eKmhs$.
\end{thm}

Theorem \ref{GPtrivial} implies that the equality
\[
\hg(\Mor_d(\pr^1, \grass))\, \hg(\mr{F}(\pr^1,n)) = \hg(\Aut(\pr^1))\, \hg(M_{0,n}(\grass,d))
\]
holds in $\nKmhs{n}$ whenever $n+3d \ge 3$. We write
\begin{equation} \label{stratoaperto}
\hg(M_{0,n}(\grass,d)) = \frac{\hg(\Mor_d(\pr^1, \grass))\, \hg(\mr{F}(\pr^1,n)) }{\hg(\Aut(\pr^1))} \,,
\end{equation}
meaning that $\hg(M_{0,n}(\grass,d))$ is the unique element of the form $\hg(X)$ such that $\hg(\Mor_d(\pr^1, \grass))\, \hg(\mr{F}(\pr^1,n)) = \hg(\Aut(\pr^1))\, \hg(X)$. 

Equation \eqref{stratoaperto} shows that in order to compute $\hg(M_{0,n}(\grass,d))$ the only missing part is the computation of $\hg(\Mor_d(\pr^1, \grass)) \in \Kmhs$. This computation is dealt with in the next section.


\section{Spaces of morphisms from $\pr^1$ to Grassmannians} \label{capmor}

The aim of this section is to show how $[\Mor_d(\pr^1, \grass)] \in \Kvar$ (and thus also $\hg(\Mor_d(\pr^1,\grass)) \in \Kmhs$) can be explicitly determined. 

Throughout the section, we work over a fixed algebraically closed field $\K$ of characteristic $0$. Schemes and morphisms between them are tacitly assumed to be over $\K$; in particular, the product is understood to be the fiber product over $\K$. Moreover, $\pr^1 = \pr^1_\K$, $V$ is a fixed $\K$-vector space of dimension $k$, and $\grass$ is the Grassmannian $G(r,V)/\K$ of $r$-dimensional quotients of $V$.

\subsection{Quot compactification and its decomposition} \label{quotcomp}

As shown in \cite{Nit05}, $\Mor_d(\pr^1, \grass)$ is naturally realized as an open subscheme of the Hilbert scheme of $\pr^1 \times \grass$, by associating to each morphism $f: \pr^1 \to \grass$ its graph $\Gamma_f \subseteq \pr^1 \times \grass$. 

A different compactification is provided by the Quot scheme
\[
\cls{Q}_d := \Quot_{V \otimes \mc{O}_{\pr^1}/ \pr^1 / \K}^{(t+1)r+d} \,,
\] 
which parametrizes equivalence classes of coherent quotients of $V \otimes \mc{O}_{\pr^1}$ of rank $r$ and degree $d$, i.e., with Hilbert polynomial equal to $(t+1)r+d \in \K[t]$. Inside $\cls{Q}_d$, $\Mor_d(\pr^1, \grass)$ is the open locus $Q_d$ corresponding to locally free quotients of $V \otimes \mc{O}_{\pr^1}$. 

The compactification $\cls{Q}_d$ was studied in \cite{Str87}. In particular, recall the following result.

\begin{thm}[{\cite[Thm.\ 2.1]{Str87}}]
The Quot scheme $\cls{Q}_d$ is an irreducible, rational, nonsingular, projective variety of dimension $k d+r(k-r)$.
\end{thm} 

The subvariety $Q_d \subseteq \cls{Q}_d$ is the open dense stratum of a certain locally closed decomposition of $\cls{Q}_d$, which we shall now study.

\begin{de}
Let $X$ be an algebraic scheme. A locally closed decomposition of $X$ is a morphism $f: Y \to X$ of algebraic schemes that satisfies the following conditions:
\begin{itemize}
\item the restriction of $f$ to each connected component of $Y$ is a locally closed immersion;
\item $f$ is bijective on closed points.
\end{itemize}
\end{de}

%
%
The starting point in constructing this decomposition of $\cls{Q}_d$ is the following observation. Since any coherent sheaf $\mc{F}$ on $\pr^1$ splits as the direct sum of its torsion subsheaf $\mc{T}(\mc{F})$ and of the locally free sheaf $\mc{F}/\mc{T}(\mc{F})$, the geometric points of $\cls{Q}_d$ parametrizing quotients that are not locally free correspond to morphisms $\pr^1 \to \grass$ of lower degree, by considering only the locally free part of the quotients. This fact suggests the possibility of finding a decomposition of $\cls{Q}_d$ such that each of its members correspond to morphisms $\pr^1 \to \grass$ of some fixed degree $\delta \le d$.

Motivated by this observation, we consider the schemes 
\[
\cls{Q}_\delta = \Quot_{V \otimes \mc{O}_{\pr^1}/ \pr^1 /\K}^{(t+1)r + \delta} 
\]
for any $\delta \in \Z$ such that $0 \le \delta \le d$. Let $[\pi_\delta: V \otimes \mc{O}_{\pr^1 \times \cls{Q}_\delta} \to \mc{F}_\delta]$ be the universal quotient on $\pr^1 \times \cls{Q}_\delta$, let $\mc{E}_\delta = \ker(\pi_\delta)$, and let $\iota_\delta: \mc{E}_\delta \to V \otimes \mc{O}_{\pr^1 \times \cls{Q}_\delta}$ be the inclusion. Then we have the short exact sequence of coherent $\mc{O}_{\pr^1 \times \cls{Q}_\delta}$-modules
\begin{equation} \label{uesQ}
0 \to \mc{E}_\delta \xto{\iota_\delta} V \otimes \mc{O}_{\pr^1 \times \cls{Q}_\delta} \xto{\pi_\delta} \mc{F}_\delta \to 0 \,,
\end{equation}
whose exactness is preserved by any base change $T \to \cls{Q}_\delta$.

Let us also introduce the relative Quot schemes
\[
\cls{R}_\delta := \Quot_{\mc{E}_\delta / \pr^1 \times \cls{Q}_\delta / \cls{Q}_\delta}^{d-\delta}
\]
over $\cls{Q}_\delta$.

\begin{prop}
The Quot scheme $\cls{R}_\delta$ is projective and smooth over $\cls{Q}_\delta$. In particular, $\cls{R}_\delta$ is nonsingular.
\end{prop}

\proof
Since the projection $\pr^1 \times \cls{Q}_\delta \to \cls{Q}_\delta$ is a projective morphism, $\cls{R}_\delta$ is a projective $\cls{Q}_\delta$-scheme (see \cite{Nit05}).

If $y \in \cls{Q}_\delta$ is a closed point, let $0 \to \mc{K} \to (\mc{E}_\delta)_y \to \mc{G}$ be the corresponding short exact sequence of $\mc{O}_{\pr^1}$-modules. Then
\[
\Ext^1(\mc{K}, \mc{G}) \cong \h^1(\pr^1, \mc{K}^\vee \otimes \mc{G}) = 0 \,,
\]
because $\mc{G}$ has $0$-dimensional support. Therefore, from \cite{Leh98} it follows that $\cls{R}_\delta$ is a smooth $\cls{Q}_\delta$-scheme.

Finally, since $\cls{Q}_\delta$ is a smooth $\K$-scheme, $\cls{R}_\delta$ is smooth over $\K$ and thus nonsingular, because $\K$ is a perfect field.
\endproof

There is a Cartesian diagram
\[
\begin{tikzcd}
\pr^1 \times \cls{R}_\delta \arrow{r}{f'_\delta} \arrow{d} & \pr^1 \times \cls{Q}_\delta \arrow{r} \arrow{d} & \pr^1 \arrow{d} \\
\cls{R}_\delta \arrow{r}{f_\delta} & \cls{Q}_\delta \arrow{r} & \spec(\K)
\end{tikzcd} ,
\]
where $f_\delta: \cls{R}_\delta \to \cls{Q}_\delta$ is the structure morphism and $f'_\delta = \id_{\pr^1} \times f_\delta$. Let
\[
0 \to (\mc{E}_\delta)_{\cls{R}_\delta} \xrightarrow{(\iota_\delta)_{\cls{R}_\delta}} V \otimes \mc{O}_{\pr^1 \times \cls{R}_\delta} \xrightarrow{(\pi_\delta)_{\cls{R}_\delta}} (\mc{F}_\delta)_{\cls{R}_\delta} \to 0
\] 
be the pullback of \eqref{uesQ} by $f'_\delta$. Moreover, let $[\rho_\delta: (\mc{E}_\delta)_{\cls{R}_\delta} \to \mc{G}_\delta]$ be the universal quotient on $\pr^1 \times \cls{R}_\delta$, and let $\mc{H}_\delta$ be the cokernel of $\ker(\rho_\delta) \into V \otimes \mc{O}_{\pr^1 \times \cls{R}_\delta}$. Then we obtain the following commutative diagram of coherent $\mc{O}_{\pr^1 \times \cls{R}_\delta}$-modules, with exact rows and columns:
\begin{equation} \label{comdiag}
\begin{tikzcd}
{} & 0 \arrow{d} & 0 \arrow{d} & 0 \arrow[dotted]{d} & \\
0 \arrow{r} & \mc{K}_\delta := \ker(\rho_\delta) \arrow{r} \arrow[equal]{d} & (\mc{E}_\delta)_{\cls{R}_\delta} \arrow{r}{\rho_\delta} \arrow{d}{(\iota_\delta)_{\cls{R}_\delta}} & \mc{G}_\delta \arrow{r} \arrow[dotted]{d} & 0 \\
0 \arrow{r} & \mc{K}_\delta \arrow{r} \arrow{d} & V \otimes \mc{O}_{\pr^1 \times \cls{R}_\delta} \arrow{r} \arrow{d}{(\pi_\delta)_{\cls{R}_\delta}} & \mc{H}_\delta \arrow{r} \arrow[dotted]{d} & 0 \\
 & 0 \arrow[dotted]{r} & (\mc{F}_\delta)_{\cls{R}_\delta} \arrow{d} \arrow[dotted]{r}{\cong} & \mc{H}_\delta/\mc{G}_\delta \arrow[dotted]{d} \arrow[dotted]{r} & 0 \\
 & & 0 & 0 &
\end{tikzcd} .
\end{equation}
The existence of the dotted arrows follows from the universal properties of kernels and cokernels. The exactness of all rows and columns is a consequence of the Four lemma and of the Nine lemma.

Since both $\mc{G}_\delta$ and $\mc{H}_\delta/\mc{G}_\delta$ are flat over $\cls{R}_\delta$, $\mc{H}_\delta$ is $\cls{R}_\delta$-flat. Thus, for any $x \in \cls{R}_\delta$ the sequence of $\mc{O}_{\pr^1_x}$-modules
\[
0 \to (\mc{G}_\delta)_x \to (\mc{H}_\delta)_x \to (\mc{H}_\delta/\mc{G}_\delta)_x \to 0
\]
is exact. Since $(\mc{H}_\delta/\mc{G}_\delta)_x \cong (\mc{F}_\delta)_{f_\delta(x)}$, it follows that the Hilbert polynomial of $(\mc{H}_\delta)_x$ is 
\[
d-\delta + (t+1)r + \delta = (t+1)r + d \in \K[t] \,.
\]
Therefore, by the universal property of $\cls{Q}_d$, there exists a unique morphism $g_\delta: \cls{R}_\delta \to \cls{Q}_d$ such that the quotient $V \otimes \mc{O}_{\pr \times \cls{R}_\delta} \to \mc{H}_\delta$ is equivalent to the pullback of $\pi_d: V \otimes \mc{O}_{\pr^1 \times \cls{Q}_d} \to \mc{F}_d$ by $\id_{\pr^1} \times g_\delta$.

\begin{rmk} \label{torslength}
The diagram \eqref{comdiag} is stable under any base change $T \to \cls{R}_\delta$, because all sheaves appearing in it are flat over $\cls{R}_\delta$. In particular, for any $x \in \cls{R}_\delta$, the fiber $(\mc{G}_\delta)_x$ is (isomorphic to) a $0$-dimensional subsheaf of $(\mc{H}_\delta)_x$, so that  there is an injective morphism from $(\mc{G}_\delta)_x$ into the torsion subsheaf $\mc{T}((\mc{H}_\delta)_x)$ of $(\mc{H}_\delta)_x$. As a consequence, the length of $\mc{T}((\mc{H}_\delta)_x)$ is at least $d-\delta$.
\end{rmk}

Let $Q_\delta$ be the open locus in $\cls{Q}_\delta$ corresponding to locally free quotients of $V \otimes \mc{O}_{\pr^1}$; $Q_\delta$ is the locus corresponding to $\Mor_\delta(\pr^1, \grass)$ inside $\cls{Q}_\delta$. Let $R_\delta$ be the preimage $f_\delta^{-1}(Q_\delta)$, with its open subscheme structure.

\begin{prop} \label{img}
For each $0 \le \delta \le d$, the set-theoretic image $g_\delta(R_\delta)$ is the constructible subset of $\cls{Q}_d$ whose closed points are the elements of
\[
\bigl\{ y \in \cls{Q}_d(\K) \mid \len\mc{T}((\mc{F}_d)_y) = d-\delta \bigr\} \,.
\]
In particular, there is a set-theoretic decomposition $\cls{Q}_d = \bigsqcup_{0 \le \delta \le d} g_\delta(R_\delta)$.
\end{prop}

\proof
For any closed point $y$ of $g_\delta(R_\delta)$, let $x \in R_\delta(\K)$ such that $g_\delta(x)=y$. 
Let us consider the pullback of \eqref{comdiag} to the fiber $\pr^1_x$:
\[
\begin{tikzcd}
{} & 0 \arrow{d} & 0 \arrow{d} & 0 \arrow{d} & \\
0 \arrow{r} & (\mc{K}_\delta)_x \arrow{r} \arrow[equal]{d} & (\mc{E}_\delta)_{f_\delta(x)} \arrow{r} \arrow{d} & (\mc{G}_\delta)_x \arrow{r} \arrow{d} & 0 \\
0 \arrow{r} & (\mc{K}_\delta)_x \arrow{r} \arrow{d} & V \otimes \mc{O}_{\pr^1} \arrow{r} \arrow{d} & (\mc{H}_\delta)_x \arrow{r} \arrow{d} & 0 \\
 & 0 \arrow{r} & (\mc{F}_\delta)_{f_\delta(x)} \arrow{d} \arrow{r}{\cong} & (\mc{H}_\delta/\mc{G}_\delta)_x \arrow{d} \arrow{r} & 0 \\
 & & 0 & 0 &
\end{tikzcd} .
\]
By Remark \ref{torslength}, $(\mc{G}_\delta)_x$ is a subsheaf of $\mc{T}((\mc{H}_\delta)_x)$. Moreover, since $x \in R_\delta$, the quotient sheaf $(\mc{H}_\delta)_x/(\mc{G}_\delta)_x \cong (\mc{H}_\delta/\mc{G}_\delta)_x \cong (\mc{F}_\delta)_{f_\delta(x)}$ is locally free. As a consequence, $(\mc{G}_\delta)_x$ actually coincides with $\mc{T}((\mc{H}_\delta)_x)$. It follows that the length of $\mc{T}((\mc{H}_\delta)_x)$ is exactly $d-\delta$. By construction of the morphism $g_\delta$, we have $(\mc{H}_\delta)_x \cong (\mc{F}_d)_y$, therefore the length of $\mc{T}((\mc{F}_d)_y)$ is $d-\delta$.

Conversely, let $y \in \cls{Q}_d$ be a closed point such that $\mc{T} = \mc{T}((\mc{F}_d)_y)$ has length $d-\delta$. Then we have a quotient
\[
V \otimes \mc{O}_{\pr^1} \to (\mc{F}_d)_y \to (\mc{F}_d)_y/\mc{T} \,,
\]
where $(\mc{F}_d)_y/\mc{T}$ is a locally free $\mc{O}_{\pr^1}$-module of rank $r$ and degree $\delta$. Let $z \in Q_\delta$ be the closed point corresponding to this quotient. There is a diagram
\[
\begin{tikzcd}
{} & 0 \arrow[dotted]{d} & 0 \arrow{d} & 0 \arrow{d} & \\
0 \arrow[dotted]{r} & \ker(\rho) \arrow[dotted]{r} \arrow[dotted]{d}{\cong} & (\mc{E}_\delta)_z \arrow[dotted]{r}{\rho} \arrow{d} & \mc{T} \arrow[dotted]{r} \arrow{d} & 0 \\
0 \arrow{r} & (\mc{E}_d)_y \arrow{r} \arrow[dotted]{d} & V \otimes \mc{O}_{\pr^1} \arrow{r} \arrow{d} & (\mc{F}_d)_y \arrow{r} \arrow{d} & 0 \\
 & 0 \arrow[dotted]{r} & (\mc{F}_\delta)_z \arrow{d} \arrow{r}{\cong} & (\mc{F}_d)_y/\mc{T} \arrow{d} \arrow{r} & 0 \\
 & & 0 & 0 &
\end{tikzcd} 
\]
which commutes and has exact rows and columns; this is a consequence of the Four lemma and of the Nine lemma. In particular, $\mc{T}$ is a $0$-dimensional quotient of $(\mc{E}_\delta)_z$ of length $d-\delta$. Therefore, $\rho: (\mc{E}_\delta)_z \to \mc{T}$ corresponds to a closed point $x \in R_\delta$ such that $f_\delta(x)=z$. Moreover, we have $g_\delta(x)=y$, hence $y \in g_\delta(R_\delta)$, as claimed.

Note that for any $y \in \cls{Q}_d$, the length of $\mc{T} = \mc{T}((\mc{F}_d)_y)$ cannot be greater than $d$. If it was, then from the quotient $V \otimes \mc{O}_{\pr^1_y} \to (\mc{F}_d)_y/\mc{T}$ we would get an invertible sheaf on $\pr^1$ with negative degree and a non-zero section. Thus, every closed point $y \in \cls{Q}_d$ lies in some $g_\delta(R_\delta)$.

By Chevalley's Theorem, each $g_\delta(R_\delta)$ is a constructible subset of $\cls{Q}_d$. Since $\K$ is algebraically closed and $\cls{Q}_d$ is of finite type over $\K$, the constructible subsets of $\cls{Q}_d$ are uniquely determined by their intersection with $\cls{Q}_d(\K)$. Therefore, from the  above discussion we deduce that $\cls{Q}_d = \bigcup_{\delta} g_\delta(R_\delta)$ and $g_\delta(R_\delta) \cap g_\epsilon(R_\epsilon) = \emptyset$ whenever $\delta \neq \epsilon$. We thus get a set-theoretic decomposition $\cls{Q}_d = \bigsqcup_{\delta} g_\delta(R_\delta)$.
\endproof

For every $0 \le \delta \le d$, define
\[
U_\delta := \cls{Q}_d \setminus \bigsqcup_{\epsilon<\delta} g_\epsilon(R_\epsilon) \,.
\]
In particular, $U_0 = \cls{Q}_d$. Let us recursively construct morphisms $j_\delta: R_\delta \to U_\delta$ as follows.

For $\delta=0$, we have $R_0 = \cls{R}_0$. Indeed, if $x \in \cls{R}_0(\K)$ then the subsheaf $(\mc{G}_0)_x \subseteq \mc{T}((\mc{H}_0)_x)$ has maximal length $d$, hence it coincides with $\mc{T}((\mc{H}_0)_x)$. Therefore, $(\mc{F}_0)_{f_0(x)} \cong (\mc{H}_0)_x/(\mc{G}_0)_x$ is locally free on $\pr^1$, i.e., $x \in R_0(\K)$. Let $j_0 := g_0: R_0 \to U_0$. Since the morphisms $g_\delta$ are projective, the image $g_0(R_0)$ is closed in $\cls{Q}_d$. Thus, $U_1 = \cls{Q}_d \setminus g_0(R_0)$ is an open subvariety of $\cls{Q}_d$.

Now, for $\delta > 0$, assume that we have already defined $j_{\delta-1}$ and proved that $U_\delta$ is an open subvariety of $\cls{Q}_d$. Since $g_\delta(R_\delta) \subseteq U_\delta$, $g_\delta|_{R_\delta}: R_\delta \to \cls{Q}_d$ factors through a morphism $j_\delta: R_\delta \to U_\delta$.

\begin{lmm} \label{cartimmer}
The commutative diagram
\[
\begin{tikzcd}
R_\delta \arrow[hookrightarrow]{r} \arrow{d}{j_\delta} & \cls{R}_\delta \arrow{d}{g_\delta} \\
U_\delta \arrow[hookrightarrow]{r} & \cls{Q}_d
\end{tikzcd}
\]
is a Cartesian diagram. As a consequence, $j_\delta: R_\delta \to U_\delta$ is projective. 
\end{lmm}

\proof
Let 
\[
\begin{tikzcd}
T \arrow{r}{u} \arrow{d}{v} & \cls{R}_\delta \arrow{d}{g_\delta} \\
U_\delta \arrow[hookrightarrow]{r} & \cls{Q}_d
\end{tikzcd}
\]
be a commutative diagram of $\K$-schemes of finite type. Let us consider the pullback of \eqref{comdiag} by $u' = \id_{\pr^1} \times u: \pr^1 \times T \to \pr^1 \times \cls{R}_\delta$. For any closed point $t \in T$, $\mc{T}((u'^\ast \mc{H}_\delta)_t)$ contains the subsheaf $(u'^\ast \mc{G}_\delta)_t$, which has length $d-\delta$. On the other hand, since $g_\delta \circ u$ factors through $U_\delta$, the length of $\mc{T}((u'^\ast \mc{H}_\delta)_t)$ cannot be greater than $d-\delta$, thus it is exactly $d-\delta$. Therefore, $\mc{T}((u'^\ast \mc{H}_\delta)_t)$ coincides with its subsheaf $(u'^\ast \mc{G}_\delta)_t$. As a consequence, the sheaf $(u'^\ast(\mc{H}_\delta/\mc{G}_\delta))_t \cong (u'^\ast\mc{H}_\delta)_t/(u'^\ast\mc{G}_\delta)_t$ is locally free.

Now, the set $\{ t \in T \mid (u'^\ast(\mc{H}_\delta/\mc{G}_\delta))_t\ \text{is locally free}\}$ is an open subset of $T$, because $u'^\ast(\mc{H}_\delta/\mc{G}_\delta)$ is $T$-flat. Since it contains all closed points of $T$, it is actually equal to $T$. Therefore, for any $t \in T$,  $((f'_\delta \circ u')^\ast\mc{F}_\delta)_t \cong (u'^\ast(\mc{H}_\delta/\mc{G}_\delta))_t$ is locally free. It follows that $u: T \to \cls{R}_\delta$ factors through a unique morphism $h: T \to R_\delta$. Since the immersion $U_\delta \into \cls{Q}_d$ is a monomorphism, we also have $j_\delta \circ h = v$.
\endproof

By Lemma \ref{cartimmer}, the image $j_\delta(R_\delta)$ is closed in $U_\delta$. Thus, the complement $U_\delta \setminus j_\delta(R_\delta)$ is an open subvariety of $\cls{Q}_d$, because $U_\delta \subseteq \cls{Q}_d$ is open. Note that we exactly have $U_{\delta+1} = U_\delta \setminus j_\delta(R_\delta)$. Thus, we can proceed recursively and define $j_\delta: R_\delta \to U_\delta$ for all $0 \le \delta \le d$.

The morphisms $j_\delta$ determine the desired decomposition of $\cls{Q}_d$.

\begin{lmm} \label{mono}
For every $0 \le \delta \le d$, $j_\delta: R_\delta \to U_\delta$ is a monomorphism.
\end{lmm}

\proof
By \cite[Prop.\ 17.2.6]{Gro67}, $j_\delta$ is a monomorphism if and only if it is radicial (equivalently, universally injective) and formally unramified. Since $j_\delta: R_\delta \to U_\delta$ is a morphism of finite type of algebraic $\K$-schemes, it suffices to prove that $j_\delta$ is injective on closed points and unramified.

Let $x_1,x_2 \in R_\delta$ be closed points such that $j_\delta(x_1) = j_\delta(x_2)$. For each $i$, pulling-back \eqref{comdiag} to $\pr^1_{x_i}$, we get the following commutative diagram with exact rows and columns:
\[
\begin{tikzcd}
{} & 0 \arrow{d} & 0 \arrow{d} & 0 \arrow{d} & \\
0 \arrow{r} & (\mc{K}_\delta)_{x_i} \arrow{r} \arrow[equal]{d} & ( (\mc{E}_\delta)_{R_\delta} )_{x_i} \arrow{r} \arrow{d} & (\mc{G}_\delta)_{x_i} \arrow{r} \arrow{d} & 0 \\
0 \arrow{r} & (\mc{K}_\delta)_{x_i} \arrow{r} \arrow{d} & V \otimes \mc{O}_{\pr^1} \arrow{r} \arrow{d} & (\mc{H}_\delta)_{x_i} \arrow{r} \arrow{d} & 0 \\
 & 0 \arrow{r} & ((\mc{F}_\delta)_{R_\delta} )_{x_i} \arrow{d} \arrow{r}{\cong} & (\mc{H}_\delta/\mc{G}_\delta)_{x_i} \arrow{d} \arrow{r} & 0 \\
 & & 0 & 0 &
\end{tikzcd} .
\]
As $j_\delta(x_1)=j_\delta(x_2)$, there is an isomorphism $(\mc{H}_\delta)_{x_1} \cong (\mc{H}_\delta)_{x_2}$ that commutes with the quotients $V \otimes \mc{O}_{\pr^1} \to (\mc{H}_\delta)_{x_i}$. Equivalently, we have $(\mc{K}_\delta)_{x_1} = (\mc{K}_\delta)_{x_2}$ as subsheaves of $\mc{O}_{\pr^1}$. The sheaves $((\mc{F}_\delta)_{R_\delta} )_{x_i}$ are locally free, because $f_\delta(x_i) \in Q_\delta$. Therefore, each $\mc{T}((\mc{H}_\delta)_{x_i})$ coincides with its subsheaf $(\mc{G}_\delta)_{x_i}$. As a consequence, we obtain that $( (\mc{E}_\delta)_{R_\delta} )_{x_1} = ( (\mc{E}_\delta)_{R_\delta} )_{x_2}$ as subsheaves of $V \otimes \mc{O}_{\pr^1}$, because both of them are the saturation of the same subsheaf $(\mc{K}_\delta)_{x_i}$. Furthermore, there is an isomorphism $(\mc{G}_\delta)_{x_1} \cong (\mc{G}_\delta)_{x_2}$ which commutes with the quotients $( (\mc{E}_\delta)_{R_\delta} )_{x_i} \to (\mc{G}_\delta)_{x_i}$, i.e., $x_1=x_2$. Thus, $j_\delta$ is injective on closed points.

Let us now prove that $j_\delta$ is unramified. Since $f_\delta|_{R_\delta}: R_\delta \to Q_\delta$ is smooth, there is a short exact sequence of locally free $\mc{O}_{R_\delta}$-modules
\begin{equation} \label{fseq}
0 \to (f_\delta|_{R_\delta})^\ast \Omega_{Q_\delta/\K} \to \Omega_{R_\delta/\K} \to \Omega_{R_\delta/Q_\delta} \to 0 \,.
\end{equation}
Moreover, we have the exact sequence of $\mc{O}_{R_\delta}$-modules
\begin{equation} \label{gseq}
g_\delta^\ast(\Omega_{U_\delta/\K}) \to \Omega_{R_\delta/\K} \to \Omega_{R_\delta/U_\delta} \to 0 \,.
\end{equation}
By \cite{Leh98}, for any closed point $x \in R_\delta$, we have
\begin{align*}
((f_\delta|_{R_\delta})^\ast \Omega_{Q_\delta/\K}) \otimes \kappa(x) &\cong \Hom\bigl( ((\mc{E}_\delta)_{R_\delta})_x, ((\mc{F}_\delta)_{R_\delta})_x \bigr)^\vee \,, \\
\Omega_{R_\delta/Q_\delta} \otimes \kappa(x) &\cong \Hom( (\mc{K}_\delta)_x, (\mc{G}_\delta)_x )^\vee \,, \\
g_\delta^\ast(\Omega_{U_\delta/\K}) \otimes \kappa(x) &\cong \Hom( (\mc{K}_\delta)_x, (\mc{H}_\delta)_x )^\vee \,.
\end{align*}
The vector space $\Hom( (\mc{G}_\delta)_x, (\mc{H}_\delta)/\mc{G}_\delta)_x )$ vanishes, because $(\mc{G}_\delta)_x$ is a torsion sheaf and $(\mc{H}_\delta)/\mc{G}_\delta)_x$ is torsion-free. Therefore, by means of the long exact Ext sequences associated to the pullback of \eqref{comdiag} to $\pr^1_x$, \eqref{fseq} implies that
\[
\Omega_{R_\delta/\K} \otimes \kappa(x) \cong \coker(\alpha)^\vee \,,
\] 
where $\alpha = - \circ (\rho_\delta)_x : \Hom((\mc{G}_\delta)_x, (\mc{H}_\delta)_x) \to \Hom(((\mc{E}_\delta)_{R_\delta})_x, (\mc{H}_\delta)_x)$. Moreover, the pullback of the morphism $g_\delta^\ast(\Omega_{U_\delta/\K}) \to \Omega_{R_\delta/\K}$ of \eqref{gseq} to $\pr^1_x$ is the dual of the canonical injection $\coker(\alpha) \into \Hom( (\mc{K}_\delta)_x, (\mc{H}_\delta)_x)$, hence it is surjective. As a consequence, we obtain the vanishing of $\Omega_{R_\delta/U_\delta} \otimes \kappa(x)$. Since $\Omega_{R_\delta/U_\delta} \otimes \kappa(x)=0$ for each closed point $x \in R_\delta$ and $R_\delta$ is of finite type over $\K$, we have $\Omega_{R_\delta/U_\delta}=0$, i.e., $j_\delta$ is unramified. 
\endproof

\begin{thm} \label{localcloseddec}
The morphism
\[
j := \coprod_{0 \le \delta \le d} \Bigl( R_\delta \xto{j_\delta} U_\delta \into \cls{Q}_d \Bigr): \coprod_{0 \le \delta \le d} R_\delta \to \cls{Q}_d
\]
is a locally closed decomposition of $\cls{Q}_d$.
\end{thm}

\proof
By Proposition \ref{img}, since $j_\delta(R_\delta) = g_\delta(R_\delta)$, $j$ is bijective on closed points. From Lemma \ref{cartimmer} and Lemma \ref{mono} it follows that $j_\delta: R_\delta \to U_\delta$ is a proper monomorphism, hence a closed immersion by \cite[Cor.\ 18.12.6]{Gro67}. Therefore, each composition $R_\delta \xto{j_\delta} U_\delta \into \cls{Q}_d$ is a locally closed immersion and $j$ is a locally closed decomposition of $\cls{Q}_d$.
\endproof

For the purposes of this work, the importance of Theorem \ref{localcloseddec} lies in the subsequent corollary, which directly follows from Proposition \ref{grringdecomp}.

\begin{cor} \label{contostrat}
The equality
\[
[ \cls{Q}_d ] = \sum_{0 \le \delta \le d} [R_\delta]
\]
holds in $\rKvar{\K}$.
\end{cor}

\subsection{The classes of the strata in $\rKvar{\K}$} \label{macosaneso}

Let us now study the relation between the classes $[R_\delta]$ and $[Q_\delta]$ in $\rKvar{\K}$.

Let $s := k-r$. By Grothendieck's theorem on vector bundles on $\pr^1$, for every $y \in Q_\delta$ there is an isomorphism
\[
(\mc{E}_\delta)_y \cong \bigoplus_{i=0}^{s} \mc{O}_{\pr^1_y}(a_i(y)) \,,
\]
where $a_1(y) \le \dots \le a_{s}(y)$ is a uniquely determined sequence of integers such that $a_1(y)+\dots+a_{s}(y)=-\delta$. Note that each $a_i(y)$ is $\le 0$, because $(\mc{E}_\delta)_y$ is a subsheaf of $V \otimes \mc{O}_{\pr^1_y}$.

Let $A_\delta = \{\mb{a} \in \Z^{s} \mid a_1 \le \dots \le a_{s} \le 0\,,\ |\mb{a}|=-\delta \}$. By \cite{Sha77}, for any $\mb{a} \in A_\delta$ the locus 
\[
Q_\delta^\mb{a} := \{ y \in Q_\delta \mid (a_1(y),\dots,a_{s}(y)) = \mb{a} \}
\]
is locally closed in $Q_\delta$, because $\mc{E}_\delta$ is flat over $Q_\delta$. Let $R_\delta^\mb{a}$ be the preimage $f_\delta^{-1}(Q_\delta^\mb{a})$, with its reduced induced scheme structure.

The main theorem of \S\ref{macosaneso} is the following.

\begin{thm} \label{loctriv}
For every $\mb{a} \in A_\delta$, let
\[
\mc{E}_\delta^\mb{a} := \bigoplus_{i=0}^{s} \mc{O}_{\pr^1}(a_i) \,.
\]
Then the restriction $f_\delta|_{R_\delta^\mb{a}}: R_\delta^\mb{a} \to Q_\delta^\mb{a}$ is a locally trivial fibration in the Zariski topology, with fiber $\Quot^{d-\delta}_{\mc{E}_\delta^\mb{a}/\pr^1/\K}$.
\end{thm}

The proof of Theorem \ref{loctriv} essentially relies on the following result.

\begin{prop} \label{pullback}
Let $T$ be a smooth scheme, and let $\mc{E}$ be a coherent locally free sheaf of rank $s$ on $\pr^1 \times T$, flat over $T$. Assume that for every $t \in T$, there is an isomorphism 
\[
\mc{E}_t \cong \bigoplus_{i=1}^s \mc{O}_{\pr^1_t}(a_i) \,,
\]
where the integers $a_1 \le \dots \le a_s$ are independent of $t$. Then, Zariski locally on $T$, $\mc{E}$ is isomorphic to the pullback of $\bigoplus_{i=1}^s \mc{O}_{\pr^1}(a_i)$ by the projection morphism $\pr^1 \times T \to \pr^1$.
\end{prop}

\proof
Without loss of generality, we may assume that $T$ is an integral scheme. Let $p: \pr^1 \times T \to \pr^1$ and $q: \pr^1 \times T \to T$ be the two projections. If $\mc{H}$ is a sheaf on $\pr^1 \times T$, we write $\mc{H}(a)$ for $\mc{H} \otimes p^\ast\mc{O}_{\pr^1}(a)$. We will proceed by induction on $s$.

If $s=1$, let us consider the coherent $\mc{O}_T$-module $\mc{F} = q_\ast( \mc{E}(-a_1))$. By the results on cohomology and base change, $\mc{F}$ is a locally free sheaf of rank $h^0(\pr^1, \mc{O}_{\pr^1})=1$, which commutes with any base change. The counit of the adjunction $q^\ast \dashv q_\ast$ gives an isomorphism of invertible sheaves $q^\ast\mc{F} \iso \mc{E}(-a_1)$. Therefore, we can cover $T$ with open subsets over which $q^\ast\mc{F}$ trivializes, so that we get an isomorphism $p^\ast\mc{O}_{\pr^1}(a_1) \iso \mc{E}$.

If $s>1$, let $\mc{F} = q_\ast(\mc{E}(-a_s))$. By Grauert's theorem on base change, $\mc{F}$ is a locally free $\mc{O}_T$-module of rank $s' = \#\{ i \mid a_i=a_s \}$. As above, the counit of the adjunction $q^\ast \dashv q_\ast$ gives a morphism of locally free sheaves $q^\ast\mc{F} \to \mc{E}(-a_s)$. On each fiber $\pr^1_t$ over $t \in T$, this is the evaluation $H^0(\pr^1_t, \mc{E}_t(-a_s)) \otimes \mc{O}_{\pr^1_t} \to \mc{E}_t(-a_s)$, which is an injection of vector bundles.
Therefore, we have an injection of locally free sheaves $(q^\ast\mc{F})(a_s) \to \mc{E}$ with locally free quotient $\mc{Q}$:
\begin{equation} \label{split}
0 \to (q^\ast\mc{F})(a_s) \to \mc{E} \to \mc{Q} \to 0 \,.
\end{equation}
On each fiber $\pr^1_t$, we have
\[
(q^\ast\mc{F})(a_s)_t \cong \bigoplus_{i \mid a_i = a_s} \mc{O}_{\pr^1_t}(a_i) \,, \qquad  \mc{Q}_t \cong \bigoplus_{i \mid a_i \neq a_s} \mc{O}_{\pr^1_t}(a_i) \,. 
\]
By the induction hypothesis, since the rank of $\mc{Q}$ is strictly less than $s$, we can cover $T$ with open subsets $U$, such that $\mc{Q}|_U$ is isomorphic to the pullback of $\bigoplus_{i \mid a_i \neq a_s} \mc{O}_{\pr^1}(a_i)$ by the projection $\pr^1 \times U \to U$. By shrinking $T$ if necessary, we may thus assume that this happens globally on $T$.

In order to prove that the sequence \eqref{split} splits (Zariski locally), consider 
\[
\Ext^1( \mc{Q}, (q^\ast\mc{F})(a_s)) \cong H^1\bigl( \pr^1 \times T, \mc{Q}^\vee \otimes (q^\ast\mc{F})(a_s) \bigr) \,.
\]
Let us apply the Leray spectral sequence $H^i(T, R^j q_\ast(-)) \Rightarrow H^{i+j}(\pr^1 \times T, -)$ to the sheaf $\mc{Q}^\vee \otimes (q^\ast\mc{F})(a_s)$. Since $\dim(\pr^1)=1$, $R^iq_\ast(\mc{Q}^\vee \otimes (q^\ast\mc{F})(a_s)) = 0$ for every $i>1$. For $i=1$, using the projection formula, we obtain
\[
R^1q_\ast(\mc{Q}^\vee \otimes (q^\ast\mc{F})(a_s)) \cong \mc{F} \otimes R^1q_\ast( \mc{Q}^\vee(a_s)) \,.
\]
The sheaf $\mc{Q}^\vee(a_s)$ is the pullback of $\bigoplus_{j \mid a_j \neq a_s} \mc{O}_{\pr^1}(a_s-a_j)$, whose $H^1$ vanishes because $a_s>a_j$, hence $R^1q_\ast( \mc{Q}^\vee(a_s)) = 0$. Finally, the last contribution to $H^1(\pr^1 \times T, \mc{Q}^\vee \otimes (q^\ast\mc{F})(a_s))$ coming from the Leray spectral sequence is 
\[
H^1\bigl( T, q_\ast(\mc{Q}^\vee \otimes (q^\ast\mc{F})(a_s)) \bigr) \cong H^1(T, \mc{F} \otimes q_\ast( \mc{Q}^\vee(a_s)) \,.
\]
If $U$ is an affine open subset of $T$, then $H^1(U, \mc{F} \otimes q_\ast( \mc{Q}^\vee(a_s)) = 0$. It follows that
\[
H^1\bigl( \pr^1 \times U, \mc{Q}^\vee \otimes (q^\ast\mc{F})(a_s) \bigr) = 0 \,,
\]
thus the sequence \eqref{split} splits over $U$. By considering an affine open cover $\{U_\lambda\}_\lambda$ such that $\mc{F}$ is trivial on each $U_\lambda$, we get the desired isomorphism.
\endproof

\proof[Proof of Theorem \ref{loctriv}]
By Proposition \ref{pullback}, there exists a Zariski open cover $\{U_\lambda\}_\lambda$ of $Q_\delta^\mb{a}$ such that $\mc{E}_\delta|_{\pr^1 \times U_\lambda} \cong p^\ast(\mc{E}_\delta^\mb{a})$ for all $\lambda$, where $p: \pr^1 \times U_\lambda \to \pr^1$ is the projection. Therefore, by the base change property of the Quot scheme, the preimage $f_\delta^{-1}(U_\lambda) \subseteq R_\delta^\mb{a}$ is naturally isomorphic to
\[
\Quot_{\mc{E}_\delta / \pr^1 \times Q_\delta^\mb{a} / Q_\delta^\mb{a}}^{d-\delta} \times_{Q_\delta^\mb{a}} U_\lambda \cong \Quot_{p^\ast(\mc{E}_\delta^\mb{a}) / \pr^1 \times U_\lambda / U_\lambda}^{d-\delta} \cong \Quot_{\mc{E}_\delta^\mb{a} / \pr^1 / \K}^{d-\delta} \times U_\lambda \,.
\]
Thus $f_\delta|_{R_\delta^\mb{a}}: R_\delta^\mb{a} \to Q_\delta^\mb{a}$ is a locally trivial fibration. 
\endproof

As a consequence of Theorem \ref{loctriv}, we get the following result.

\begin{cor} \label{contotriv}
For each $0 \le \delta \le d$ and each $\mb{a} \in A_\delta$, the equality
\[
[R_\delta] = \sum_{\mb{a} \in A_\delta} [ Q_\delta^\mb{a} ] \Bigl[ \Quot_{\mc{E}_\delta^\mb{a} / \pr^1 / \K}^{d-\delta} \Bigr]
\]
holds in $\rKvar{\K}$.
\end{cor}

\proof
The subschemes $R_\delta^\mb{a}$, for $\mb{a} \in A_\delta$, give a locally closed decomposition of $R_\delta$. Therefore, we have
\[
[R_\delta] = \sum_{\mb{a} \in A_\delta} [ R_\delta^\mb{a} ] \,.
\]
By Theorem \ref{loctriv}, each $f_\delta|_{R_\delta^\mb{a}}: R_\delta^\mb{a} \to Q_\delta^\mb{a}$ is a locally trivial fibration in the Zariski topology, with fiber $\Quot_{\mc{E}_\delta^\mb{a} / \pr^1 / \K}^{d-\delta}$. Thus, from Proposition \ref{solitarobatriviale} it follows that
\[
[ R_\delta^\mb{a} ] = [ Q_\delta^\mb{a} ] \Bigl[ \Quot_{\mc{E}_\delta^\mb{a} / \pr^1 / \K}^{d-\delta} \Bigr] \,.
\]
Combining the two equalities, we obtain the desired result.
\endproof

The problem of computing the class of $\Quot_{\mc{E}_\delta^\mb{a} / \pr^1 / \K}^{d-\delta}$ in $\rKvar{\K}$ was solved in \cite{BFP19}, where it was proved that it is independent of $\mb{a} \in A_\delta$.

\begin{prop} \label{contoquot}
For any integer $0 \le \delta \le d$, the equality
\[
\Bigl[ \Quot_{\mc{E}_\delta^\mb{a} / \pr^1 / \K}^{d-\delta} \Bigr] = \sum_{\substack{ \mb{m} \in (\Z_{\ge 0})^{s} \\ |\mb{m}|=d-\delta }} \frac{ ( 1-\maff^{m_1 +1} ) \cdots ( 1-\maff^{m_{s} +1} ) }{ (1-\maff)^{s} } \maff^{d_\mb{m}}
\]
holds in $\rKvar{\K}$, where $d_\mb{m} = \sum_{i=1}^{s} (i-1)m_i$. 
\end{prop}

\proof
Since
\[
[ \pr^m ] = \sum_{j=0}^{m} \maff^j = \frac{1-\maff^{m+1}}{1-\maff} \,,
\]
the result follows from the equality
\[
\Bigl[ \Quot_{\mc{E}_\delta^\mb{a} / \pr^1 / \K}^{d-\delta} \Bigr] = \sum_{\substack{ \mb{m} \in (\Z_{\ge 0})^{s} \\ |\mb{m}|=d-\delta }} [\pr^{m_1}] \cdots [\pr^{m_{s}} ]\, \maff^{d_\mb{m}} \,.
\]
of \cite[Prop. 4.5]{BFP19}.
\endproof

The class appearing in Proposition \ref{contoquot}, which depends only on $s=k-r$ and $d-\delta$, will be denoted by $\Omega_{d-\delta}$.

\subsection{The class of $\Mor_d(\pr^1, G(r,V))$ in $\rKvar{\K}$} \label{scarmor}

By means of the previous results, we can now describe a method for computing the class $[Q_d] = [\Mor_d(\pr^1, \grass)] \in \rKvar{\K}$.

\begin{thm} \label{morkvar}
The equality
\[
\sum_{d=0}^\infty [Q_d] q^d = \Biggl( \sum_{d=0}^\infty [ \cls{Q}_d ] q^d \Biggr) \Biggl( \sum_{d=0}^\infty \Omega_d q^d \Biggr)^{-1}
\]
holds in $\rKvar{\K}\ser$. 
\end{thm}

\proof
Applying Corollary \ref{contotriv} and Proposition \ref{contoquot}, we get
\[
[R_\delta] = \sum_{\mb{a} \in A_\delta} [ Q_\delta^\mb{a} ] \Bigl[ \Quot_{\mc{E}_\delta^\mb{a} / \pr^1 / \K}^{d-\delta} \Bigr] = \Omega_{d-\delta} \Biggl( \sum_{\mb{a} \in A_\delta} [ Q_\delta^\mb{a} ] \Biggr) = \Omega_{d-\delta} [Q_\delta] \,.
\]
By Corollary \ref{contostrat}, we thus have
\[
[ \cls{Q}_d ] = \sum_{\delta=0}^d [R_\delta] = \sum_{\delta=0}^d \Omega_{d-\delta} [Q_\delta] \,,
\]
therefore
\begin{equation} \label{giggino}
\sum_{d=0}^\infty [ \cls{Q}_d ] q^d = \Biggl( \sum_{d=0}^\infty [Q_d] q^d \Biggr) \Biggl( \sum_{d=0}^\infty \Omega_d q^d \Biggr) \,.
\end{equation}
Since $\Omega_0 = 1$, the last power series in \eqref{giggino} is invertible, thus \eqref{giggino} implies the theorem.
\endproof

Theorem \ref{morkvar} provides a recursive formula for $[Q_d]$ in terms of $[\cls{Q}_\delta]$ and $\Omega_\delta$ for $\delta \le d$:
\begin{equation} \label{rec}
\mho_0 := 1 \,, \quad \mho_j := -\sum_{i=0}^{j-1} \mho_i \Omega_{j-i} \,, \quad [Q_d] = \sum_{j=0}^{d} \mho_j [ \cls{Q}_{d-j} ] \,.
\end{equation}
The elements $\Omega_\delta$ have already been determined in Proposition \ref{contoquot}. Thus, in order for \eqref{rec} to be effective for computing $[Q_d]$, we only need to calculate $[\cls{Q}_\delta]$.

By \cite[\S 3]{Str87}, $\cls{Q}_\delta$ has a Białynicki-Birula decomposition determined by a torus action. Since each subvariety in this locally closed decomposition is isomorphic to some affine space $\aff^i$, we have
\begin{equation} \label{bigquotmotive}
[\cls{Q}_\delta] = \sum_{i=0}^{\dim(\cls{Q}_\delta)} m_{\delta,i} [\aff^i] = \sum_{i=0}^{k\delta + r(k - r)} m_{\delta,i} \maff^i \,,
\end{equation}
where $m_{\delta,i}$ is the number of $i$-dimensional cells. 

A formula for the numbers $m_{\delta,i}$ is provided in \cite{Str87}. Let $G \subseteq \Z^{s} \times \Z^{s+1} \times \Z^{s}$ be the subset whose elements are triples $(\mb{a}=(a_1,\dots,a_{s})$, $\mb{b}=(b_0,\dots,b_{s})$, $\mb{c}=(c_1,\dots,c_{s}))$ such that
\begin{gather*}
b_0 = 0 \le a_1 \le b_1 \le a_2 \le \dots \le b_{s-1} \le a_{s} \le \delta = b_{s} \,, \\
0 \le c_1 \le \dots \le c_{s} \le r \,.
\end{gather*}
Then the following result holds.

\begin{prop}[{\cite[Thm.\ 5.4]{Str87}}] \label{strommecomp}
For all $0 \le i \le k \delta+r (k - r)$, $m_{\delta,i}$ is equal to the number of elements $(\mb{a},\mb{b},\mb{c}) \in G$ such that the equality
\[
\sum_{j=1}^{s} (a_j + c_j (1 + b_j - b_{j-1}) ) = i 
\]
holds.
\end{prop}

By means of this proposition, \eqref{bigquotmotive} becomes an explicit formula for the class $[\cls{Q}_\delta] \in \rKvar{\K}$. This completely solves the problem of calculating $[Q_d]$: after computing $[\cls{Q}_\delta]$ for all $\delta \le d$ via \eqref{bigquotmotive}, $[Q_d]$ can be determined via \eqref{rec} and Proposition \ref{contoquot}.

\section{The Betti numbers of $\cls{M}_{0,n}(G(r,V),d)$} \label{esempiacaso}

The previous results yield an algorithm to compute the HG-characteristic of $\cls{M}_{0,n}(\grass,d)$, which in turn determines its Hodge and Betti numbers. Let us outline how this algorithm works.
\begin{enumerate}[label=(\roman*)]
\item Apply the results of \S\ref{scarmor} to compute $[Q_\delta]$ and $\hg(Q_\delta)$ for all $\delta \le d$.
\item Calculate $\hg(M_{0,m}(\grass,\delta))$ for all $\delta \le d$ and $m \le n+d-\delta$, by means of \eqref{stratoaperto}.
\item Use the recursive procedure explained in \S\ref{sectonrecrel} to get $\hg(\cls{M}_{0,n}(\grass,d))$.
\end{enumerate}

\begin{rmk}
As a consequence of the above procedure, the cohomology group $H^i(\cls{M}_{0,n}(\grass,d))$ vanishes for $i$ odd, whereas its class in $\Kmhs$ is a multiple of $\caff^{i/2}$ for $i$ even. This agrees with the results of \cite{Opr06}.
\end{rmk}

In the remainder of the section, we examine some examples, in order to show the effectiveness of our algorithm in practice.

\subsection{The HG-characteristic of $\cls{M}_{0,0}(G(2,4), 2)$} \label{G(2,4)2}

Let us consider the Grassmannian $\grass = G(2,4)$, whose class in $\Kvar$ is
\[
[G(2,4)] = (\maff^2 + 1) (\maff^2 + \maff + 1) \,.
\]
First, we compute $\Omega_i$ for all $0 \le i \le 2$ via Proposition \ref{contoquot}:
\[
\Omega_0 = 1 \,, \qquad
\Omega_1 = (\maff + 1)^2 \,, \qquad
\Omega_2 = \maff^4 + 2\maff^3 + 4\maff^2 + 2\maff + 1 \,.
\]
Then we also have the following equalities:
\[
\mho_0 = 1 \,, \qquad
\mho_1 = -(\maff + 1)^2 \,, \qquad
\mho_2 = 2\maff (\maff^2 + \maff + 1) \,.
\]
The classes $[\cls{Q}_\delta]$ for all $0 \le \delta \le 2$ are determined by Proposition \ref{strommecomp}:
\begin{align*}
[\cls{Q}_0] &= (\maff^2 + 1) (\maff^2 + \maff + 1) \,, \\
[\cls{Q}_1] &= (\maff + 1)^2 (\maff^2 + 1) (\maff^2 - \maff + 1) (\maff^2 + \maff + 1) \,, \\
[\cls{Q}_2] &= (\maff^2 + 1) (\maff^2 + \maff + 1) (\maff^4 + 1) (\maff^4 + \maff^3 + \maff^2 + \maff + 1) \,.
\end{align*}
Combining the previous equalities, we get
\begin{align*}
[Q_0] &= (\maff^2 + 1) (\maff^2 + \maff + 1) \,, \\
[Q_1] &= \maff (\maff - 1) (\maff + 1)^2 (\maff^2 + 1) (\maff^2 + \maff + 1) \,,
\\
[Q_2] &= \maff^3 (\maff - 1) (\maff + 1) (\maff^2 + 1) (\maff^2 + \maff + 1) (\maff^3 + \maff^2 + \maff - 1) \,.
\end{align*}

The knowledge of $\hg(Q_d)$ allows us to calculate $\hg (M_{0,n}(\grass, d) )$. By \eqref{stratoaperto} and Corollary \ref{configserrecar}, for all $d>0$ we have
\[
\hg(M_{0,0}(\grass,d)) = \frac{ \hg(Q_d) }{ \caff^3 - \caff } \qquad
\text{and} \qquad
\hg(M_{0,1}(\grass,d)) =  \frac{ \hg(Q_d) }{ \caff^2 - \caff } s_1 \,,
\]
therefore
\begin{align*}
\hg(M_{0,m}(\grass,0)) &= \hg(\cls{M}_{0,m}(\grass,0)) = 0 \qquad \text{for all}\ m \le 2 \,, \\
\hg(M_{0,0}(\grass,1)) &= \hg(\cls{M}_{0,0}(\grass,1)) = (\caff + 1) (\caff^2 + 1) (\caff^2 + \caff + 1) \,, \\
\hg( M_{0,1}(\grass,1) ) &= \hg( \cls{M}_{0,1}(\grass,1) ) = (\caff + 1)^2 (\caff^2 + 1) (\caff^2 + \caff + 1) s_1 \,, \\
\hg(M_{0,0}(\grass,2)) &= \caff^2 (\caff^2 + 1) (\caff^2 + \caff + 1) (\caff^3 + \caff^2 + \caff - 1) 
\,. 
\end{align*}
By \eqref{recfibre}, we also have
\[
\hg(\cls{\Phi}_{0,0}) = 0 \qquad
\text{and} \qquad
\hg(\cls{\Phi}_{0,1}) = \frac{D(\hg(M_{0,1}(\grass,1)))}{\hg(\grass)} = (\caff + 1)^2 \,.
\]

Now, we can conclude our computation by means of Corollary \ref{veryimpcoroll}. Since $p_i \circ -$ is an algebra homomorphism, Corollary \ref{veryimpcoroll} and the equations $s_1 = p_1$ and $2 s_2 = p_2 + p_1^2$ imply that
\begin{multline} \label{summ00Y1}
\hg( \cls{M}_{0,0}(\grass,2) ) = \hg( M_{0,0}(\grass,2) ) + \frac{ \hg( M_{0,1}(\grass,1) ) }{ s_1 } \hg( \cls{\Phi}_{0,1} ) \\
+ \hg(\grass) \frac{ p_2 \circ \hg( \cls{\Phi}_{0,1} ) - \hg( \cls{\Phi}_{0,1} )^2 }{ 2 } \,.
\end{multline}
Thus, the above calculations yield
\[
\hg( \cls{M}_{0,0}(\grass,2) ) = \caff^9 + 3\caff^8 + 7\caff^7 + 11\caff^6 + 14\caff^5 + 14\caff^4 + 11\caff^3 + 7\caff^2 + 3\caff + 1 \,.
\]
As a consequence, the $E$-polynomial $E(\cls{M}_{0,0}(\grass,2); t,u)$ is
\[
t^9u^9 + 3t^8u^8 + 7t^7u^7 + 11t^6u^6 + 14t^5u^5 + 14t^4u^4 + 11t^3u^3 + 7t^2u^2 + 3tu + 1 \,,
\]
in agreement with \cite[Thm.\ 3.1]{Lop14}.

\subsection{The HG-characteristic of $\cls{M}_{0,1}(G(2,4), 2)$}

By \eqref{stratoaperto} and Corollary \ref{configserrecar}, we have
\begin{align*}
\hg( M_{0,3}(\grass,0) ) &= (\caff^2 + 1) (\caff^2 + \caff + 1) s_3 \,, \\
\hg(M_{0,2}(\grass,1)) &= \caff (\caff + 1) (\caff^2 + 1) (\caff^2 + \caff + 1) \bigl( (\caff - 1) s_2 + s_1^2 \bigr) \,, \\
\hg( M_{0,1}(\grass,2) ) &= \caff^2 (\caff + 1) (\caff^2 + 1) (\caff^2 + \caff + 1) (\caff^3 + \caff^2 + \caff - 1) s_1 \,.
\end{align*}
As a consequence, we get
\[
\hg( \cls{\Phi}_{1,1} ) = \frac{ D( \hg( M_{0,2}(\grass,1) ) ) }{ \hg(\grass) } + \frac{ D( \hg( M_{0,3}(\grass,0) ) ) }{ \hg(\grass)\, s_2 } \hg( \cls{\Phi}_{0,1} ) s_1 = (\caff + 1)^3 s_1 \,.
\]

Using these equalities and those of \S\ref{G(2,4)2}, we can now compute $\hg( \cls{M}_{0,1}(\grass,2) )$. Corollary \ref{veryimpcoroll} implies that $\hg( \cls{M}_{0,1}(\grass,2) )$ is equal to
\begin{multline*}
\hg( M_{0,1}(\grass,2) ) + \frac{ \hg( M_{0,1}(\grass,1) ) }{ s_1 } \hg( \cls{\Phi}_{1,1} )+ \{ \hg( M_{0,2}(\grass,1) ) \}_{s_2}\, \hg( \cls{\Phi}_{0,1} ) s_1 \\
+ 2 \{ \hg( M_{0,2}(\grass,1) ) \}_{s_1^2}\, \hg( \cls{\Phi}_{0,1} ) s_1 + \frac{ \hg( M_{0,3}(\grass,0) ) }{ 2 s_3 } \hg(\cls{\Phi}_{0,1})^2 s_1 \\
+ \frac{ \hg( M_{0,3}(\grass,0) ) }{ 2 s_3 } \bigl( p_2 \circ \hg(\cls{\Phi}_{0,1}) \bigr) s_1 - \hg(\grass)\, \hg(\cls{\Phi}_{0,1})\, \hg( \cls{\Phi}_{1,1} ) \,,
\end{multline*}
where $\{ \hg( M_{0,2}(\grass,1) ) \}_{s_2}$ (resp.\ $\{ \hg( M_{0,2}(\grass,1) ) \}_{s_1^2}$) is the coefficient of $s_2$ (resp.\ $s_1^2$) in $\hg( M_{0,2}(\grass,1) )$. Therefore, the HG-characteristic of $\cls{M}_{0,1}(\grass,2)$ is equal to
\[
( \caff^{10} + 4\caff^9 + 12\caff^8 + 22\caff^7 + 33\caff^6 + 36\caff^5 + 33\caff^4 + 22\caff^3 + 12\caff^2 + 4\caff + 1 ) s_1 \,.
\]
In particular, its $E$-polynomial $E(\cls{M}_{0,1}(\grass,2); t,u)$ is
\begin{multline*}
t^{10}u^{10} + 4t^9u^9 + 12t^8u^8 + 22t^7u^7 + 33t^6u^6 + 36t^5u^5 + 33t^4u^4 + 22t^3u^3 \\ 
+ 12t^2u^2 + 4tu + 1 \,.
\end{multline*}

\bibliographystyle{amsalpha}
\bibliography{bib}

\providecommand{\bysame}{\leavevmode\hbox to3em{\hrulefill}\thinspace}
\providecommand{\MR}{\relax\ifhmode\unskip\space\fi MR }
\providecommand{\MRhref}[2]{%
  \href{http://www.ams.org/mathscinet-getitem?mr=#1}{#2}
}
\providecommand{\href}[2]{#2}
\begin{thebibliography}{ACG11}

\bibitem[ACG11]{ACG11}
E.~Arbarello, M.~Cornalba, and {P.\,A.} Griffiths, \emph{Geometry of algebraic
  curves. {V}ol.\ {II}}, Grund.\ Math.\ Wiss., vol. 268, Springer, Heidelberg,
  2011.

\bibitem[Bag19]{Bag19}
M.~Bagnarol, \emph{On the cohomology of moduli spaces of stable maps to
  {G}rassmannians}, Ph.D. thesis, SISSA, Trieste, 2019,
  http://hdl.handle.net/20.500.11767/103198.

\bibitem[BFP19]{BFP19}
M.~Bagnarol, B.~Fantechi, and F.~Perroni, \emph{On the motive of {Q}uot schemes
  of zero-dimensional quotients on a curve}, arXiv e-prints (2019),
  arXiv:1907.00826.

\bibitem[Bit04]{Bit04}
F.~Bittner, \emph{The universal {E}uler characteristic for varieties of
  characteristic zero}, Comp.\ Math. \textbf{140} (2004), no.~4, 1011–1032.

\bibitem[BLL97]{BLL97}
F.~Bergeron, G.~Labelle, and P.~Leroux, \emph{Combinatorial species and
  tree-like structures}, Encycl.\ Math.\ Appl., Cambridge Univ.\ Press,
  Cambridge, UK, 1997.

\bibitem[BM96]{BM96}
K.~Behrend and Yu. Manin, \emph{Stacks of stable maps and {G}romov-{W}itten
  invariants}, Duke Math.\ J. \textbf{85} (1996), no.~1, 1--60.

\bibitem[BO03]{BO03}
K.~Behrend and A.~O'Halloran, \emph{On the cohomology of stable map spaces},
  Invent.\ Math. \textbf{154} (2003), no.~2, 385--450.

\bibitem[Bor18]{Bor18}
L.~Borisov, \emph{The class of the affine line is a zero divisor in the
  grothendieck ring}, J.\ Alg.\ Geom. \textbf{27} (2018), 203--209.

\bibitem[DK86]{DK86}
{V.\,I.} Danilov and {A.\,G.} Khovanski\u{\i}, \emph{Newton polyhedra and an
  algorithm for calculating {H}odge-{D}eligne numbers}, Izv. Akad. Nauk SSSR
  Ser. Mat. \textbf{50} (1986), no.~5, 925--945.

\bibitem[FP97]{FP97}
W.~Fulton and R.~Pandharipande, \emph{Notes on stable maps and quantum
  cohomology}, Algebraic {G}eometry {S}anta {C}ruz 1995, Proc.\ Sympos.\ Pure
  Math., vol.~62, Amer.\ Math.\ Soc., Providence, RI, 1997, pp.~{45--96}.

\bibitem[Get95]{Get95conf}
E.~Getzler, \emph{Mixed {H}odge structures of configuration spaces}, arXiv
  e-prints (1995), alg--geom/9510018.

\bibitem[GP06]{GP06}
E.~Getzler and R.~Pandharipande, \emph{The {B}etti numbers of
  {$\overline{\mathcal{M}}_{0,n}(r,d)$}}, J.\ Alg.\ Geom. \textbf{15} (2006),
  no.~4, 709--732.

\bibitem[Gro67]{Gro67}
A.~Grothendieck, \emph{{\'E}l\'ements de g\'eom\'etrie alg\'ebrique: {IV}.
  \'{E}tude locale des sch\'emas et des morphismes de sch\'emas, {Q}uatri\`eme
  partie}, Publ.\ Math.\ IH\'ES \textbf{32} (1967), 5--361.

\bibitem[Joy81]{Joy81}
A.~Joyal, \emph{Une théorie combinatoire des séries formelles}, Adv.\ Math.
  \textbf{42} (1981), no.~1, 1--82.

\bibitem[Kel05]{Kel05}
{G.\,M.} Kelly, \emph{On the operads of {{J}.{P}.\ {M}ay}}, Repr.\ Theory
  Appl.\ Categ. \textbf{13} (2005), 1--13.

\bibitem[Kon95]{Kon95}
M.~Kontsevich, \emph{Enumeration of rational curves via torus action}, The
  Moduli Space of Curves ({T}exel {I}sland, 1994), Progr.\ Math., vol. 129,
  Birkh\"{a}user Boston, Cambridge, MA, 1995, pp.~335--368.

\bibitem[Leh98]{Leh98}
M.~Lehn, \emph{On the cotangent sheaf of {Q}uot-schemes}, Internat.\ J.\ Math.
  \textbf{9} (1998), no.~4, 513--522.

\bibitem[{Ló}14]{Lop14}
A.~{López Martín}, \emph{{P}oincaré polynomials of stable map spaces to
  {G}rassmannians}, Rend.\ Sem.\ Mat.\ Univ.\ Padova \textbf{131} (2014),
  193--208.

\bibitem[Mac95]{Mac95}
{I.\,G.} Macdonald, \emph{Symmetric functions and {H}all polynomials}, 2nd ed.,
  Oxford Math.\ Mon., The Clarendon Press, Oxford Univ.\ Press, New York, 1995.

\bibitem[MM07]{Mu07}
A.~Musta{\c{t}}{\v{a}} and M.A. Musta{\c{t}}{\v{a}}, \emph{Intermediate moduli
  spaces of stable maps}, Invent.\ Math. \textbf{167} (2007), no.~1, 47--90.

\bibitem[Nit05]{Nit05}
N.~Nitsure, \emph{Construction of {H}ilbert and {Q}uot schemes}, Fundamental
  Algebraic Geometry, Math.\ Surveys Monogr., vol. 123, Amer. Math. Soc.,
  Providence, RI, 2005, pp.~105--137.

\bibitem[Opr06]{Opr06}
D.~Oprea, \emph{The tautological rings of the moduli spaces of stable maps to
  flag varieties}, J.\ Alg.\ Geom. \textbf{15} (2006), no.~4, 623--655.

\bibitem[PS08]{PS08}
C.~Peters and {J.\,H.\,M.} Steenbrink, \emph{Mixed {H}odge structures}, Erg.\
  der Math., vol.~52, Springer, Berlin, 2008.

\bibitem[Sha77]{Sha77}
{S.\,S.} Shatz, \emph{The decomposition and specialization of algebraic
  families of vector bundles}, Comp.\ Math. \textbf{35} (1977), no.~2,
  163--187.

\bibitem[Str87]{Str87}
{S.\,A.} Str{\o}mme, \emph{On parametrized rational curves in {G}rassmann
  varieties}, Space Curves ({R}occa di {P}apa, 1985), Lect.\ Notes Math., vol.
  1266, Springer, Berlin, Heidelberg, 1987, pp.~251--272.

\end{thebibliography}

\Addresses

\end{document}